\documentclass{cmip-l}
\usepackage[all]{xy}
\usepackage{amssymb}
\usepackage{amsmath}
\swapnumbers
\newtheorem{thm}[subsection]{Theorem}
\newtheorem{prop}[subsection]{Proposition}
\newtheorem{lemma}[subsection]{Lemma}

\newenvironment{proo}{\begin{trivlist} \item{\emph{Proof.}}}
  {\hfill $\square$ \end{trivlist}}

\def\d{\succ}
\def\g{\prec}

\def\r{\vdash }

\def\DD{\Delta}

\def\ss{\sigma}

\def\t{\otimes}

\def\Vtn{V^{\otimes n}}

\def\gg{\mathfrak{g}}

\newcommand{\Vt}[1]{V^{\otimes #1}}

\def\AA{{\mathcal{A}}}

\def\pt{\textrm{-}}

\def\HH{{\mathcal{H}}}
\def\CCC{{\mathcal{C}}}

\def\PP{{\mathcal{P}}}

\def\QQQ{{\mathcal{Q}}}
\def\RRR{{\mathcal{R}}}
\def\SSS{{\mathcal{S}}}
\def\XXX{{\mathcal{X}}}
\def\YYY{{\mathcal{Y}}}
\def\ZZZ{{\mathcal{Z}}}

\def\Prim{\mathrm{Prim\, }}
\def\Id{\mathrm{Id }}

\def\epi{\twoheadrightarrow}
\def\mono{\rightarrowtail}

\def\Bi{\mathop{B_{\infty}}}

\def\Ker{\mathop{\rm Ker }}

\def\dim{\mathop{\rm dim }}
\def\id{\mathrm{ id }}

\def\PTi{PT_{\infty}}

\def\DDo{\overline {\Delta}}
\def\HHo{\overline {\HH}}
\def\To{\overline {T}}
\def\So{\overline {S}}
\def\Ao{\overline {A}}

\def\alg{\textrm{-alg}}

\def\KK{\mathbb{K}}

\def\sms{{\bf S}-modules }

\def\arbreA{\vcenter{\xymatrix@R=3pt@C=3pt{
&& \\
&*{}\ar@{-}[ur] \ar@{-}[ul] \ar@{-}[d]     &\\
&&
}}}

\def\arbreAv{\vcenter{\xymatrix@R=3pt@C=3pt{
&v& \\
&*{}\ar@{-}[ur] \ar@{-}[ul] \ar@{-}[d]     &\\
&&
}}}

\def\arbreBA{\vcenter{\xymatrix@R=2pt@C=2pt{
&&&&\\
&&&*{}\ar@{-}[ul] & \\
&&*{}\ar@{-}[uurr] \ar@{-}[uull] \ar@{-}[d]     &&\\
&&&&
}}}

\def\arbreBAv{\vcenter{\xymatrix@R=2pt@C=2pt{
&v_{1}&&v_{2}&\\
&&&*{}\ar@{-}[ul] & \\
&&*{}\ar@{-}[uurr] \ar@{-}[uull] \ar@{-}[d]     &&\\
&&&&
}}}

\def\arbreAB{\vcenter{\xymatrix@R=2pt@C=2pt{
&&&&\\
&*{}\ar@{-}[ur] &&& \\
&&*{}\ar@{-}[uurr] \ar@{-}[uull] \ar@{-}[d]     &&\\
&&&&
}}}

\def\arbreABv{\vcenter{\xymatrix@R=2pt@C=2pt{
&v_{1}&&v_{2}&\\
&*{}\ar@{-}[ur] &&& \\
&&*{}\ar@{-}[uurr] \ar@{-}[uull] \ar@{-}[d]     &&\\
&&&&
}}}

\def\arbreBB{\vcenter{\xymatrix@R=2pt@C=2pt{
&&*{}&&\\
&&&& \\
&&*{}\ar@{-}[uurr] \ar@{-}[uull] \ar@{-}[d] \ar@{-}[uu]     &&\\
&&&&
}}}

\def\arbreABC{\vcenter{\xymatrix@R=1pt@C=1pt{
&&&&&&\\
&*{}\ar@{-}[ur] &&&&& \\
&&*{}\ar@{-}[uurr] &&&&\\
&&&*{}\ar@{-}[uuurrr] \ar@{-}[uuulll] \ar@{-}[d] &&&\\
&&&&&&
}}}

\def\arbreBAC{\vcenter{\xymatrix@R=1pt@C=1pt{
&&&&&&\\
&&&*{}\ar@{-}[ul] &&& \\
&&*{}\ar@{-}[uurr] &&&&\\
&&&*{}\ar@{-}[uuurrr] \ar@{-}[uuulll] \ar@{-}[d] &&&\\
&&&&&&
}}}

\def\arbreACA{\vcenter{\xymatrix@R=1pt@C=1pt{
&&&&&&\\
&*{}\ar@{-}[ur] &&&&*{}\ar@{-}[ul] & \\
&&&&&&\\
&&&*{}\ar@{-}[uuurrr] \ar@{-}[uuulll] \ar@{-}[d] &&&\\
&&&&&&
}}}

\def\arbreACAv{\vcenter{\xymatrix@R=1pt@C=1pt{
&v_{1}&&v_{2}&&v_{3}&\\
&*{}\ar@{-}[ur] &&&&*{}\ar@{-}[ul] & \\
&&&&&&\\
&&&*{}\ar@{-}[uuurrr] \ar@{-}[uuulll] \ar@{-}[d] &&&\\
&&&&&&
}}}

\def\arbreACB{\vcenter{\xymatrix@R=1pt@C=1pt{
&&&&&&\\
&*{}\ar@{-}[ur] &&&&& \\
&&&&*{}\ar@{-}[uull] &&\\
&&&*{}\ar@{-}[uuurrr] \ar@{-}[uuulll] \ar@{-}[d] &&&\\
&&&&&&
}}}

\def\arbreBCA{\vcenter{\xymatrix@R=1pt@C=1pt{
&&&&&&\\
&&&&&*{}\ar@{-}[ul] & \\
&&*{}\ar@{-}[uurr] &&&&\\
&&&*{}\ar@{-}[uuurrr] \ar@{-}[uuulll] \ar@{-}[d] &&&\\
&&&&&&
}}}

\def\arbreCAB{\vcenter{\xymatrix@R=1pt@C=1pt{
&&&&&&\\
&&&*{}\ar@{-}[ur] &&& \\
&&&&*{}\ar@{-}[uull] &&\\
&&&*{}\ar@{-}[uuurrr] \ar@{-}[uuulll] \ar@{-}[d] &&&\\
&&&&&&
}}}

\def\arbreCBA{\vcenter{\xymatrix@R=1pt@C=1pt{
&&&&&&\\
&&&&&*{}\ar@{-}[ul] & \\
&&&&*{}\ar@{-}[uull] &&\\
&&&*{}\ar@{-}[uuurrr] \ar@{-}[uuulll] \ar@{-}[d] &&&\\
&&&&&&
}}}

\def\unreducedtrees{
\begin{center}
\begin{picture}(350,70)

\put(10,00){\scriptsize $\bullet$}

\put(30,00){\scriptsize$\bullet$}
\put(30,20){\scriptsize $\bullet$}
\put(32,00){\line(0,1){20}}

\put(50,00){\scriptsize $\bullet$}
\put(50,20){\scriptsize $\bullet$}
\put(50,40){\scriptsize $\bullet $}
\put(52,00){\line(0,1){20}}
\put(52,20){\line(0,1){20}}

\put(90,00){\scriptsize $\bullet$}
\put(70,20){\scriptsize $\bullet$}
\put(110,20){\scriptsize $\bullet$}
\put(92,00){\line(-1,1){20}}
\put(90,00){\line(1,1){20}}

\put(130,00){\scriptsize $\bullet$}
\put(130,20){\scriptsize $\bullet$}
\put(130,40){\scriptsize $\bullet $}
\put(130,60){\scriptsize $\bullet $}
\put(132,00){\line(0,1){20}}
\put(132,20){\line(0,1){20}}
\put(132,40){\line(0,1){20}}

\put(160,00){\scriptsize$\bullet$}

\put(160,20){\scriptsize $\bullet$}
\put(162,00){\line(0,1){20}}
\put(140,40){\scriptsize $\bullet$}
\put(180,40){\scriptsize $\bullet$}
\put(162,20){\line(-1,1){20}}
\put(160,20){\line(1,1){20}}

\put(210,00){\scriptsize$\bullet$}
\put(190,20){\scriptsize $\bullet$}
\put(214,00){\line(-1,1){20}}
\put(230,40){\scriptsize $\bullet$}
\put(230,20){\scriptsize $\bullet$}
\put(232,20){\line(0,1){20}}
\put(210,00){\line(1,1){20}}

\put(270,00){\scriptsize$\bullet$}
\put(250,20){\scriptsize $\bullet$}
\put(274,00){\line(-1,1){20}}
\put(250,40){\scriptsize $\bullet$}
\put(290,20){\scriptsize $\bullet$}
\put(252,20){\line(0,1){20}}
\put(270,00){\line(1,1){20}}

\put(330,00){\scriptsize $\bullet$}
\put(310,20){\scriptsize $\bullet$}
\put(330,20){\scriptsize $\bullet$}
\put(332,00){\line(0,1){20}}
\put(350,20){\scriptsize $\bullet$}
\put(332,00){\line(-1,1){20}}
\put(330,00){\line(1,1){20}}
\end{picture}
\end{center}
}

\def\unreducedtreesdeux{
\begin{center}
\begin{picture}(100,45)

\put(10,00){\scriptsize$\bullet$}
\put(-10,20){\scriptsize $\bullet$}
\put(14,00){\line(-1,1){20}}
\put(30,40){\scriptsize $\bullet$}
\put(30,20){\scriptsize $\bullet$}
\put(32,20){\line(0,1){20}}
\put(10,00){\line(1,1){20}}

\put(70,00){\scriptsize$\bullet$}
\put(50,20){\scriptsize $\bullet$}
\put(74,00){\line(-1,1){20}}
\put(50,40){\scriptsize $\bullet$}
\put(90,20){\scriptsize $\bullet$}
\put(52,20){\line(0,1){20}}
\put(70,00){\line(1,1){20}}

\end{picture}
\end{center}
}

\def\labeledunreducedtrees{
\begin{center}
\begin{picture}(350,80)

\put(10,00){\scriptsize $\bullet$}
\put(5,00){\scriptsize 1 }

\put(30,00){\scriptsize$\bullet$}
\put(25,00){\scriptsize 1}
\put(30,20){\scriptsize $\bullet$}
\put(25,20){\scriptsize 2}
\put(32,00){\line(0,1){20}}

\put(50,00){\scriptsize $\bullet$}
\put(45,00){\scriptsize 1}
\put(50,20){\scriptsize $\bullet$}
\put(45,20){\scriptsize 2}
\put(50,40){\scriptsize $\bullet $}
\put(45,40){\scriptsize 3}
\put(52,00){\line(0,1){20}}
\put(52,20){\line(0,1){20}}

\put(90,00){\scriptsize $\bullet$}
\put(85,00){\scriptsize 1}
\put(70,20){\scriptsize $\bullet$}
\put(70,25){\scriptsize 2}
\put(110,20){\scriptsize $\bullet$}
\put(110,25){\scriptsize 3}
\put(92,00){\line(-1,1){20}}
\put(90,00){\line(1,1){20}}

\put(130,00){\scriptsize $\bullet$}
\put(125,00){\scriptsize 1}
\put(130,20){\scriptsize $\bullet$}
\put(125,20){\scriptsize 2}
\put(130,40){\scriptsize $\bullet $}
\put(125,40){\scriptsize 3}
\put(130,60){\scriptsize $\bullet $}
\put(125,60){\scriptsize 4}
\put(132,00){\line(0,1){20}}
\put(132,20){\line(0,1){20}}
\put(132,40){\line(0,1){20}}

\put(160,00){\scriptsize$\bullet$}
\put(155,00){\scriptsize 1}
\put(160,20){\scriptsize $\bullet$}
\put(155,20){\scriptsize 2}
\put(162,00){\line(0,1){20}}
\put(140,40){\scriptsize $\bullet$}
\put(140,45){\scriptsize 3}
\put(180,40){\scriptsize $\bullet$}
\put(180,45){\scriptsize 4}
\put(162,20){\line(-1,1){20}}
\put(160,20){\line(1,1){20}}

\put(210,00){\scriptsize$\bullet$}
\put(200,00){\scriptsize 1}
\put(190,20){\scriptsize $\bullet$}
\put(185,20){\scriptsize 2}
\put(214,00){\line(-1,1){20}}
\put(190,40){\scriptsize $\bullet$}
\put(190,45){\scriptsize 3}
\put(230,20){\scriptsize $\bullet$}
\put(230,25){\scriptsize 4}
\put(192,20){\line(0,1){20}}
\put(210,00){\line(1,1){20}}

\put(270,00){\scriptsize$\bullet$}
\put(260,00){\scriptsize 1}
\put(250,20){\scriptsize $\bullet$}
\put(245,20){\scriptsize 2}
\put(274,00){\line(-1,1){20}}
\put(290,40){\scriptsize $\bullet$}
\put(290,45){\scriptsize 4}
\put(290,20){\scriptsize $\bullet$}
\put(285,25){\scriptsize 3}
\put(292,20){\line(0,1){20}}
\put(270,00){\line(1,1){20}}

\put(330,00){\scriptsize $\bullet$}
\put(320,00){\scriptsize 1}
\put(310,20){\scriptsize $\bullet$}
\put(310,25){\scriptsize 2}
\put(330,20){\scriptsize $\bullet$}
\put(330,25){\scriptsize 3}
\put(332,00){\line(0,1){20}}
\put(350,20){\scriptsize $\bullet$}
\put(350,25){\scriptsize 4}
\put(332,00){\line(-1,1){20}}
\put(330,00){\line(1,1){20}}
\end{picture}
\end{center}
}

\begin{document}

\author[J.-L. Loday, M. Ronco]{Jean-Louis Loday, Mar\'ia Ronco}
\address{JLL: Institut de Recherche Math\'ematique Avanc\'ee,
    CNRS et Universit\'e de Strasbourg,
    7 rue R. Descartes,
   67084 Strasbourg Cedex, France}
\email{loday@math.u-strasbg.fr}
\address{MOR: Departamento de Matematicas,
  Facultad de Ciencias-Universidad de Valparaiso,
      Avda. Gran Bretana 1091,
         Valparaiso, Chile}
\email{maria.ronco@uv.cl}

\title{Combinatorial Hopf algebras}


\thanks{This work has been partially supported by  the FONDECYT Project 1085004, the ``Ecos-Sud'' project E01C06, an ``Agence Nationale pour la Recherche'' (ANR) contract.}
\date{\today}

\dedicatory{Ce travail est d\'edi\'e \`a Alain Connes en t\' emoignage d'estime et d'une longue amiti\'e}

\begin{abstract} 

We define a  ``combinatorial Hopf algebra'' as a Hopf algebra which is free (or cofree) and equipped with a given isomorphism to the free algebra over the indecomposables (resp.\ the cofree coalgebra over the primitives). The choice of such an isomorphism implies the existence of a finer algebraic structure on the Hopf algebra and on the indecomposables (resp.\  the primitives). For instance a cofree-cocommutative right-sided combinatorial Hopf algebra is completely determined by its primitive part which is a pre-Lie algebra.  The key example is the Connes-Kreimer Hopf algebra. The study of all these combinatorial Hopf algebra types gives rise to several good triples of operads. It involves the operads: dendriform, pre-Lie, brace, and variations of them.

\end{abstract}

\maketitle

\section*{Introduction} \label{s:int}

Many recent papers are devoted to some infinite dimensional Hopf algebras called collectively ``combinatorial Hopf algebras''. The Connes-Kreimer Hopf algebra is one of them \cite{CK}. Among other examples we find:  the Fa\`a di Bruno algebra, variations of the symmetric functions: $Sym, NSym, FSym, QSym, FQSym,  PQSym,$ (Aguiar-Sottile \cite{AguiarSottile, AguiarBergeronSottile}, Bergeron-Hohlweg \cite{BergeronHohlweg}, Chapoton-Livernet \cite{ChapotonLivernet07}, Livernet  \cite{Livernet06}, Malvenuto-Reutenauer  \cite{MR}, Hivert-Novelli-Thibon  \cite{HNT08}, Palacios-Ronco \cite{PalaciosRonco}), examples related to knot theory (Turaev  \cite{Turaev}), to quantum field theory (Brouder  \cite{Brouder}, Figueroa-Gracia-Bondia  \cite{FGB}, Brouder-Frabetti-Krattenthaler \cite{BrouderFrabetti, BFK}, van Suijlekom \cite{Suijlekom}), to foliations (Connes-Moscovici \cite{ConnesMoscovici}), to $K$-theory (Gangl-Goncharov-Levin  \cite{GGL}, Loday-Ronco  \cite{LRHopf}, Lam-Pylyavskyy  \cite{LamPylyavskyy}).
The aim of this paper is to make precise the meaning of combinatorial Hopf algebras and to unravel their fine algebraic structure.

Here is a typical example of a combinatorial Hopf algebra. The Connes-Kreimer algebra $\HH_{CK}$ is a free-commutative Hopf algebra over some families of trees $t$. The coalgebra structure is of the form
$$\DD(t)=\sum_{c}P_{c}(t) \t R_{c}(t),$$
where $P_{c}(t) $ is a polynomial of trees and $R_{c}(t) $ is a tree (the sum is over admissible cuts $c$, cf.\ \ref{freecomrsCHA}). In this description there are two key points. First, we observe that, not only 
$\HH_{CK}$ is free as a commutative algebra, but we are given a precise isomorphism with the polynomial algebra over the indecomposables (spanned by the trees). This is why we call it ``combinatorial'' Hopf algebra. Second, the coproduct has a special form since its second component is linear (a tree instead of a polynomial of trees). This is called the ``right-sided condition''. A priori the space of indecomposables is a Lie coalgebra. But it is known, in this case, that it has a finer structure: it is a \emph{pre-Lie} coalgebra. Our claim is that the existence of this pre-Lie structure is not special to the Connes-Kreimer algebra, but is a general fact for combinatorial Hopf algebras satisfying the right-sided condition.

This case (free-commutative Hopf algebra) can be dualized (in the graded sense)  to give rise to the following statement: a Hopf algebra structure on the cofree-cocommutative coalgebra $S^c(V)$, which satisfies some ``right-sided condition'' (cf.\ \ref{conditionr-s}), induces a pre-Lie algebra structure on $V$ (cf.\ Theorem \ref{thm:cocomCHApreLie}). The Lie product on the primitive part $V$ is the anti-symmetrization of this pre-Lie product.

 We will see that there are several contexts for combinatorial Hopf algebras, depending on the following choice of options: free or cofree, associative or commutative, general or right-sided.
Since the free and cofree cases are dual to one another under graded dualization we will only study one of them in details, namely the cofree case. Some generalizations are evoked in the last section.

In the cofree-coassociative general context a combinatorial Hopf algebra is a Hopf algebra structure on the tensor coalgebra $T^c(R)$ for a certain vector space $R$. We show that, for  such a combinatorial Hopf algebra, the primitive part $R$ is a \emph{multibrace algebra} ($MB$-algebra) and that any multibrace algebra $R$ gives rise to a Hopf structure on $T^c(R)$. A multibrace algebra is determined by $(p+q)$-ary operations $M_{pq}$ satisfying some relations analogous to the relations satisfied by the brace algebras (see \ref{Rijk}). Moreover $T^c(R)$ inherits a finer algebraic structure: it is a \emph{dipterous algebra} (see \ref{dipterousstructure}). So, in the cofree-coassociative context, the classification of the combinatorial Hopf algebras involves the operad $\PP= MB$ governing the primitive part and the operad $\AA= Dipt$ governing the algebra structure. We show that the combinatorial Hopf algebra associated to a free multibrace algebra is in fact a free dipterous algebra:
$$Dipt(V) \cong T^c(MB(V)).$$
 This is part of a more general result which says that there is a ``good triple of operads'' (in the sense of  \cite{JLLgbto}):
$$(Ass, Dipt, MB),$$
where $Ass$ is the operad of associative algebras.

Our aim is to handle similarly the three other cases (cofree-coassociative right-sided, cofree-cocommutative general, cofree-cocommutative right-sided) and to determine the operad $\PP$ and the operad $\AA$ when possible. We show that these two operads are strongly related since there is a good triple of operads
$$(\CCC, \AA, \PP)$$
where $\CCC= Ass$ or $Com$ (governing commutative algebras) depending on the context we are working in.
We show that the operads $\PP$ and $\AA$ are as follows:

\bigskip
\centerline{
\begin{tabular}{| c ||  c | c  c c | c |  }
\hline
		&  $\CCC$ &&& $\AA$ & $\PP$ \\
\hline
                 &           & $Dipt$&& &  \\
   general  &  $As$ & &or&        & $MB$  \\
                 &           & &&$2as$ &  \\
\hline
                 &           & $Dend$&& &  \\
right-sided &  $As$ & &or&        & $Brace$  \\
                  &           & &&$\YYY:=2as/\sim$ &  \\
\hline
\hline
   general  &  $Com$ & &&$ComAs$        & $SMB$  \\
\hline
right-sided &  $Com$ & &&$\XXX:= ?$     & $SBrace=PreLie$  \\
           \hline
\end{tabular}
}

\bigskip
In the last case we conjecture the existence of an operad structure on $Com\circ preLie$, that we denote by $\XXX$.
Some of these operads are known, the other ones are variations of known types:

\begin{itemize}
\item $MB$ : multibrace algebra, also called ${\bf B_{\infty}}$-algebra in \cite{LRstr}, non-differential $\Bi$ in \cite{GerstenhaberVoronov}, Hirsch algebra in \cite{Kadeishvili}, $LR$-algebra in \cite{MerkulovVallette}, see \ref{multibrace},
\item $Brace$ : brace algebra \cite{GerstenhaberVoronov, GJ}, see \ref{brace},

\item $SMB$ : symmetric multibrace algebra, a symmetric variation of $MB$, see \ref{SMB},

\item $PreLie$ : pre-Lie algebra, also called Vinberg algebra, or right-symmetric algebra, see \ref{rsccCHA},

\item $Dipt$ : dipterous algebra \cite{LRnote}, see \ref{dipterousstructure},

\item $Dend$ : dendriform algebra  \cite{JLLdig}, see \ref{dend},

\item $2as$ : 2-associative algebra (two associative products) \cite{LRstr}, see \ref{2as},

\item $\YYY:=2as/\{M_{pq}=0\  |\  p\geq 2\}$, see \ref{ss:YYY},

\item $ComAs$ : commutative-associative algebra (a commutative product and an associative product), see \ref{ComAs}.

\end{itemize}

The dimension of the space of $n$-ary operations is given in the following tableau:

\bigskip

\centerline{
\begin{tabular}{| c ||  c | c | c |  }
\hline
		&  $\CCC$ & $\AA$ & $\PP$ \\
\hline\hline
   general, $As$  &  $n!$ &  $2 C_{n-1}\times n!$         & $C_{n}\times n!$  \\
\hline
right-sided, $As$ &  $n!$ &$c_{n}\times n!$    & $c_{n-1}\times n!$  \\
\hline
\hline
   general, $Com$  &  $1$ & $d_{n}$        & $f_{n}$  \\
\hline
right-sided, $Com$ &  $1$ & $(n+1)^{n-1}$     & $n^{n-1}$  \\
           \hline
\end{tabular}
}

\bigskip
\noindent where $c_{n}$ is the Catalan number, $C_{n}$ is the super-Catalan number, $d_{n}$ is the number of  ``labeled series-parallel posets with $n$ points'' and  $f_{n}$ is the number of  ``connected labeled series-parallel posets with $n$ nodes'' (cf.\ \cite{Stanley74, StanleyVol2} problem 5.39 and \ref{operadSMB}):

\medskip

\centerline{
\begin{tabular}{| c ||  c | c | c | c | c | c| }
\hline
	$n$	&  1& 2& 3 & 4   & 5& $\cdots$ \\
\hline
\hline
$c_{n}$ & 1 & 2 & 5 & 14 & 42& $\cdots$ \\
\hline
$C_{n}$ & 1 & 3 & 11& 45 & 197& $\cdots$\\
\hline
$\dim ComAs(n)=d_{n}$ & 1 & 3 & 19 & 195 & 2791& $\cdots$ \\
\hline
$\dim SMB(n)=f_{n}$ & 1 & 2 & 12 & 122 & 1740& $\cdots$ \\
           \hline
\end{tabular}
}

\bigskip

We observe that some of these operads are rather simple since they are binary and quadratic: $PreLie, Dipt, Dend, 2as, ComAs$. The others are more complicated since they involve $n$-ary generating operations for any $n$: $MB,  Brace, SMB$, or nonquadratic relations: $\YYY$. 

Summarizing our results, we have shown the existence of the following good triples of operads:

$$\begin{array}{ccc}
(As, Dipt, MB),&(As, 2as, MB),&(Com, ComAs, SMB),\\
(As, Dend, Brace),&
(As, \YYY, Brace),&
(Com, \XXX, PreLie),
\end{array}$$
except for the last one which is conjectural.
\bigskip

In most examples the space of primitives $R$ comes endowed with a basis of ``combinatorial objects'' like permutations, trees, graphs, tableaux. Therefore the Hopf algebra is made of polynomials on these combinatorial objects. For instance in every context the free algebra $\PP(V)$ can be described by means of trees. Changing the basis of $R$ does not change the $\PP$-algebra structure of $R$. However, if the coalgebra isomorphism $\HH\cong T^c(R)$ (resp.\ $S^c(R)$) is modified, then the $\PP$-algebra structure of $R$ is modified accordingly. For instance, in the cofree-cocommutative right-sided context, the structure of pre-Lie algebra of $R$ is modified (but not its structure of Lie algebra). In \ref{MalvenutoReutenauer} we make explicit two different combinatorial Hopf algebra structures on the Malvenuto-Reutenauer Hopf algebra. This importance of the basis had been envisioned by Joni and Rota in their seminal paper \cite{JoniRota} where they say:\emph{ ``It must be stressed that the coalgebras of combinatorics come equipped with a distinguished basis, and many an interesting combinatorial problem can be formulated algebraically as that of transforming this basis into another basis with more desirable properties. Thus, a mere structure theory of coalgebras---or Hopf algebras---will hardly be sufficient for combinatorial purposes.''}

\bigskip

The plan of the paper is as follows. In the first section we recall the basic notion of Hopf algebra and the results on triples of operads that are going to be used in the proofs.
The other four sections are devoted to the following four cases:

\begin{itemize}
\item section 2: cofree-coassociative CHA and $MB$-algebras,
\item section 3: cofree-coassociative right-sided CHA  and $Brace$-algebras,
\item section 4: cofree-cocommutative CHA  and $SMB$-algebras,
\item section 5: cofree-cocommutative  right-sided CHA  and $PreLie$-algebras.
\end{itemize}

The plan of each of these four sections is as follows:

\begin{itemize}
\item definition of the CHAs involved,
\item the algebraic structure of the primitives, the operad $\PP$, the equivalence,
\item the algebraic structure of the Hopf algebra, the operad $\AA$,
\item comparison of $\PP$ and $\AA$, the good triple $(\CCC, \AA, \PP)$, where $\CCC=As$ or $Com$,
\item combinatorial description of the free algebras $\PP(V)$ and $\AA(V)$ when available,
\item the dual case (free-associative or free-commutative),
\item variation (in the associative context)
\end{itemize}

The first case (cofree-coassociative general CHA) is treated in details. When the proofs in the other cases are analogous they are, most of the time, omitted. In the last section we list several examples from the literature.
\medskip

Finally, let us say a word about the terminology. In the literature the term ``combinatorial Hopf algebras'' is used to call Hopf algebras based on combinatorial objects (cf.\ for instance \cite{HivertNovelliThibon}), without a precise definition about what is a combinatorial object. In all cases these Hopf algebras are free (or cofree) and the combinatorial objects provide a basis, whence an isomorphism as required in our definition. The only exception is the use of combinatorial Hopf algebra in the paper \cite{AguiarBergeronSottile} where an extra piece of information is required (namely a character map).

\medskip

\noindent{\bf Acknowledgement.} Thanks to Emily Burgunder, Vladimir Dotsenko, Ralf Holtkamp, Muriel Livernet and Bruno Vallette for useful comments on earlier versions of this paper. Thanks also to the referee whose comments were very helpful in clarifying a few points.

\bigskip

\noindent {\bf Notation.} In this paper $\KK$ is a field and all vector spaces 
are over $\KK$. Its unit is denoted by 1. 
The vector space spanned by the elements of a set $X$ is denoted $\KK [X]$. 
The tensor product  of vector spaces over $\KK$
is denoted by $\t$. The tensor product of $n$ copies of the space $V$ is
denoted $\Vtn$. For $v_i\in V$ the element $v_1\t \cdots \t v_n$ of 
$\Vtn$ is denoted by the concatenation of the elements: $v_1 \ldots v_n$.
The \emph{tensor module} over $V$ is the direct sum
$$T(V) := \KK 1\oplus V\oplus V^{\t 2}\oplus \cdots \oplus V^{\t n} \oplus \cdots $$
and the \emph{reduced tensor module} is $\To(V) := T(V)/\KK 1$.
The \emph{symmetric module} over $V$ is the direct sum
$$S(V) := \KK 1\oplus V\oplus S^2(V) \oplus \cdots \oplus S^n(V) \oplus \cdots $$
and the \emph{reduced symmetric module} is $\bar S(V) := S(V)/\KK 1$,
where $S^n(V)= (V^{\t n})_{S_{n}}$ is the quotient of $V^{\t n}$ by the action of the symmetric group. We still denote by $v_1 \ldots v_n$ the image in $S^n(V)$ of $v_1 \ldots v_n\in V^{\t n}$.
If $V$ is generated by $x_1,\ldots, x_n$, then $S(V)$ (resp. $T(V)$) can be identified with the polynomials (resp. noncommutative polynomials) in $n$ variables.

If the set $X$ is a basis of $V=\KK[X]$, then we write $T(X)$ (resp. $S(X)$) in place of $T(V)$ (resp. $S(V)$).

 \section{Prerequisites: Hopf algebras and operads} 
 
 \subsection{Hopf algebras} Recall that a \emph{bialgebra} is a vector space $\HH$ equipped with an associative and unital algebra structure $(\HH , *, u)$ and a coassociative counital coalgebra structure $(\HH , \DD, \epsilon)$, which satisfy the \emph{Hopf compatibility relation}, that is $\Delta: \HH \to \HH \t \HH$ and  $\epsilon: \HH \to \KK $ are morphisms of associative unital algebras.
 
 It is helpful to introduce the augmentation ideal $\HHo:= \Ker \epsilon$ and the reduced comultiplication 
 $\DDo: \HHo\to \HHo\t \HHo$ defined by the formula $\DD(x) = x\t 1 + 1\t x + \DDo(x)$. 
 The iteration of the map $\DDo$ gives rise to $n$-ary cooperations $\DDo^{n} : \HHo \to \HHo^{\t n}$. The filtration of $\HHo$ is defined by:
$$F_r\HHo := \{x\in \HHo \ \vert \ \DDo^{n}(x)=0 \textrm{ for any }   n>r \}.$$
By definition $\HH$ is said to be \emph{connected} (or \emph{conilpotent})  if $\HHo= \bigcup_{r\geq 1}F_{r}\HHo$.
Observe that the first piece of the filtration is the space of primitives of the bialgebra: $F_1\HHo = \Prim \HH$. Since any conilpotent bialgebra can be equipped with an antipode, there is an equivalence between conilpotent Hopf algebras and conilpotent bialgebras. 

\subsection{Cofree bialgebras} We say that a conilpotent bialgebra is \emph{cofree-coasso\-ciative} if, as a conilpotent coalgebra, it is a cofree object. Denoting by $R$ the space of primitives, it means that there exists an isomorphism of coalgebras $\HH \cong T^c{R}$. Recall that the cofree coalgebra $T^c{R}$ is the tensor module as a graded vector space, and that the coproduct is the \emph{deconcatenation}:
$$\DD(x_{1}\cdots x_{n})= \sum_{i=0}^{i=n}x_{1}\cdots x_{i}\t x_{i+1}\cdots x_{n}.$$
As a result the following universal property holds: any conilpotent coalgebra homomorphism $C \to T^c(R)$ is completely determined by the composite $C \to T^c(R) \to R$. On the other hand, any linear map $\varphi : C\to R$ which maps $1$ to $0$ determines a unique coalgebra homomorphism $C\to T^c(R)$. The component in $R^{\t n}$ of the image of $c\in C$ is $\sum \varphi(c_{(1)})\t \cdots \t \varphi(c_{(n)})$, where $\DD^{n}(c) = \sum c_{(1)}\t \cdots \t c_{(n)}$.

Similar statements hold in the cocommutative case with $S^c(R)$ in place of $T^c(R)$.

\subsection{Operads}\label{ops} For an operad $\PP$ the free $\PP$-algebra over the vector space $V$ is denoted $\PP(V)$. It is of the form $\PP(V)= \bigoplus_{n}\PP(n)\t _{S_{n}}V^{\t n}$ where $\PP(n)$ is the space of $n$-ary operations considered as a right module. The left module structure of $V^{\t n}$ is given by $\ssÊ\dot (v_{1}\cdots v_{n}) = v_{\ss^{-1}(1)}\cdots v_{\ss^{-1}(n)}$ . If $\PP(n)$ is free as a representation of $S_{n}$, then we write it $\PP(n)=\PP_{n}\t \KK[S_{n}]$ and $\PP_{n}$ is called the space of nonsymmetric $n$-ary operations. In this case we obtain $\PP(V) =  \bigoplus_{n}\PP_{n}\t V^{\t n}$.

The generating series of the operad $\PP$ is 
$$f^{\PP}(t) := \sum_{n\geq 1}\dim \frac{\PP(n)}{n!} t^n.$$
The operads governing the associative algebras, commutative algebras, Lie algebras, pre-Lie algebras are denoted $As,\ Com,\ Lie,\ PreLie$ respectively. For more on operads, see for instance \cite {LV}.

\subsection{Triple of operads}\label{triples} Recall briefly from \cite{JLLgbto} what it means for $(\CCC, \AA, \PP)$ to be a good triple of operads. First, there is supposed to be a well-defined notion of $\CCC^c\pt \AA$-bialgebra, $\CCC$ governing the coalgebra structure and $\AA$ governing the algebra structure. The operations and the cooperations are supposed to be related by \emph{compatibility relations}. In our cases these compatibility relations are always distributive, that is the composite of an operation with a cooperation can be written as the composite of a cooperation followed by an operation (or an algebraic sum of them). Then,

-- the primitive part of any $\CCC^c\pt \AA$-bialgebra is a $\PP$-algebra,

\medskip

-- there is a pair of adjoint functors ($U\dashv F$): $\AA\alg {{U \atop \longleftarrow} \atop {\longrightarrow\atop F}} \PP\alg$,

\medskip

-- the following structure theorem holds: 

\noindent for any  $\CCC^c\pt \AA$-bialgebra $\HH$ the following are equivalent:

\medskip

(a) $\HH$ is connected,
 
 (b)  $\HH$ is isomorphic to $U(\Prim \HH)$,
 
(c)  $\HH$ is cofree over its primitive part, i.e.\ isomorphic to $\CCC^c(\Prim \HH)$.

\medskip

An important consequence of this theorem is an isomorphism of vector spaces  $\AA(V)\cong \CCC^c(\PP(V))$, which is functorial in $V$.

In \cite{JLLgbto} we found a criterion which ensures that the notion of $\CCC^c\pt\AA$-bialgebra gives rise to a good triple of operads $(\CCC, \AA, \PP)$, where the operad $\PP:=\Prim_{\CCC}\AA$ is made of the primitive operations. There are three conditions:

\noindent (H0) The compatibility relations are distributive.

\noindent (H1) The free $\AA$-algebra $\AA(V)$ is naturally a $\CCC^c\pt\AA$-bialgebra.

\noindent (H2epi) The $\CCC$-coalgebra map $\AA(V)\to \CCC^c(V)$ (deduced from H1) admits a splitting.

\medskip 

The isomorphism $\AA(V)\cong \CCC^c(\PP(V))$, which is a consequence of these hypotheses, depends on the choice of the splitting. It implies $f^{\AA}(t)=f^{\CCC}(f^{\PP}(t))$.

In the cases that we are looking at in this paper $\CCC$ is either the operad $As$ or the operad $Com$. The condition (H0) is immediate to verify by direct inspection. So the first task (once $\AA$ is discovered) consists in verifying the hypotheses (H1) and (H2epi). The second task consists in unravelling the operad $\PP$ by providing a presentation by generators and relations.

\section{(Co)free-(co)associative CHA}\label{s:cofreecoas}

In this section we study the cofree-coassociative CHAs. Some of the results of this section were announced in \cite{LRnote}, where the notion of dipterous algebras was introduced for the first time. Examples will be given in section \ref{Exs}. They include  the Malvenuto-Reutenauer algebra (\ref{MalvenutoReutenauer}), the incidence Hopf algebras (\ref{incidence}) and the Hopf algebra on a tensor module (\ref{tensormodule}). 

\subsection{Definition} A \emph{cofree-coassociative combinatorial Hopf algebra} $\HH$ (cofree-coassociative CHA for short) is a cofree bialgebra with a prescribed isomorphism with $T^c(R)$ (here $R:=\Prim \HH$). Equivalently, it is a vector space $R$ with a product $*$ on $T^c(R)$ which makes $(T^c(R), *, \DD)$ ($\DD=$ deconcatenation) into a Hopf algebra.

It is important to notice that the isomorphism $\HH\cong T^c(R)$ is part of the structure. Any change of this isomorphism changes the product $*$. Such an isomorphism is completely determined by a splitting of the inclusion $R=\Prim \HH \rightarrowtail \HH$ (analogue of a Lie idempotent).

A morphism of cofree-coassociative CHAs is a linear map $\varphi : R \to R'$ whose extension  $T^c\varphi : T^c (R) \to T^c (R')$ is a morphism of Hopf algebras. 

\subsection{Multibraces} Let  $\HH=T^c(R)$ be a cofree-coassociative CHA. Since it is cofree, the multiplication $*: T^c(R)\t T^c(R)\to T^c(R)$ is completely determined by its projection onto $R$. So, there are well-defined maps $M_{pq}: R^{\t p}\t R^{\t q} \to R$, called \emph{multibrace operations}, given by the following composite:

$$\xymatrix{
T^c(R)\t T^c(R)\ar^-{*}[r] & T^c(R)\ar^{proj}[d] \\
 R^{\t p}\t R^{\t q} \ar^{inc}[u]\ar^-{M_{pq}}[r] & R
}$$ 

It is sometimes helpful to write
$$\{x_{1}, \ldots , x_{p}; y_{1}, \ldots , y_{q}\}:= M_{pq}(x_{1} \ldots  x_{p} ; y_{1} \ldots , y_{q})\ .$$
  One also finds the notation $\{x_{1}, \ldots , x_{p}\}\{ y_{1}, \ldots , y_{q}\}$ in the literature \cite{GerstenhaberVoronov, Akman}.
  
  \begin{prop}\label{Rijk} Let $M_{pq}$ be $(p+q)$-ary operations on $R$ ($p\geq 0, q\geq 0$), and let $*$ be the unique binary operation on $T^c(R)$ which is a coalgebra morphism and whose projection onto $R$ coincides with the $M_{pq}$'s. The following assertions are equivalent:
  
  (a) the operation $*$ is associative and unital,
  
  (b) the operations $M_{pq}$ satisfy the following relations $\RRR$:
  $$M_{00}=0, M_{01}= \id = M_{10}, M_{0q}= 0= M_{p0}, \textrm { for } p>1, q>1,$$
  and, for any integers $i,j,k$ greater than or equal to $1$,
  
  $ (\RRR_{ijk}):$
  $$ \sum _{1\leq l \leq i+j} M_{lk}\circ (M_{i_{1} j_{1}}\ldots M_{i_{l} j_{l}} ;\id^{\t k}) =
   \sum _{1\leq m \leq j+k} M_{im}\circ (\id^{\t i} ;M_{j_{1} k_{1}}\ldots M_{j_{m} k_{m}})$$
where the left sum is extended to all sets of indices $i_{1}, \ldots , i_{l}; j_{1}, \ldots , j_{l}$ such that $i_{1}+ \cdots + i_{l}= i ; j_{1}+ \cdots + j_{l}= j$, and  the right sum is extended to all sets of indices $ j_{1}, \ldots , j_{m}; k_{1}, \ldots , k_{m}$ such that $j_{1}+ \cdots + j_{m}= j ; k_{1}+ \cdots + k_{m}= k$.  \end{prop}

In the formula the concatenation of operations $M_{pq}M_{rs}$ has the following meaning:
$$\displaylines{
M_{pq}M_{rs}(x_{1}\ldots x_{p+r} ; y_{1}\ldots y_{q+s}) =\hfill \cr
\hfill M_{pq}(x_{1}\ldots x_{p} ; y_{1}\ldots y_{q}) M_{rs}(x_{p+1}\ldots x_{p+r} ; y_{q+1}\ldots y_{q+s}).
}$$

  \begin{proo} See  \cite {LRnote} and Proposition 1.6 of \cite {LRstr}.
  \end{proo}
  
  \subsection{Example}\label{FormulaforM11} The first nontrivial relation is $(\RRR_{111})$ which reads
  $$M_{21}(uv+vu,w) + M_{11}(M_{11}(u,v), w) =  M_{11}(u, M_{11}(v,w))+ M_{12}(u, vw+wv)$$
  for any $u,v,w\in R$.
  
  \subsection{Multibrace algebra}\label{multibrace} By definition a \emph{ multibrace algebra} (or $MB$-algebra) is a vector space $R$ equipped with multibrace operations $M_{pq}$, that is $p+q$-ary operations $M_{pq}:R^{\t p}\t R^{\t q} \to R$, which satisfy all the conditions $\RRR$ of Proposition \ref{Rijk}.
  
  This notion was called ${\bf B}_{\infty}$-algebra in \cite{LRstr} and first appeared in the differential graded framework as $\Bi$-algebra in \cite{GerstenhaberVoronov, GJ}. It appears in \cite{MerkulovVallette} as $LR$-algebras.
  
  \begin{thm}\label{coAsCombHopfMB} There is an equivalence of categories between the cofree-coasso\-ciative combinatorial Hopf algebras and the multibrace algebras.
  \end{thm}
  \begin{proo} Taking the primitives gives a functor $\Prim: coAs\pt CHA\to MB\alg$. In the other direction the functor $MB\alg\to coAs\pt CHA$, $R\mapsto T^c(R)$,  is a consequence of Proposition \ref{Rijk}. These two functors are inverse to each other.
  \end{proo}
  
  \subsection{Dipterous structures}\label{dipterousstructure} 
  A \emph{dipterous algebra} (cf.\ \cite{LRnote}) is a vector space $A$ equipped with two binary operations $*$ and $\d$ which satisfy the two relations
 $$(x*y)* z = x* (y* z), \quad \textrm{and} \quad (x*y)\d z = x\d (y\d z)$$
  for any $x,y, z \in A$. In most examples the associative algebras that we are working with are unital and augmented. Therefore it is helpful to introduce the notion of \emph{unital dipterous algebra}. It is an augmented algebra $A=\KK \ 1 \oplus \Ao$ such that $\Ao$ is a dipterous algebra, i.e.\ equipped with a right operation. We ask that the following relations hold for any $a\in \Ao$:
  $$ a\d 1 = 0, \quad \textrm{and} \quad 1\d a = a.$$
  Observe that $1\d 1$ is not defined, but 1 is a unit for $*$.
  
  If $A$ and $B$ are two dipterous algebras, then there is a way to construct a dipterous structure on the tensor product $A\t B$ as follows. The associative product is as usual:
  $$(a\t b) * (a'\t b') = (a * a') \t (b * b').$$
  The right product is given by
$$(a\t b) \d (a'\t b') = (a * a') \t (b \d b').$$
  This formula makes sense provided that we do not have $b=1$ and $b'=1$ simultaneously. In that case we put
  $$(a\t 1)\d (a'\t 1) = (a\d a')\t 1\ .$$
  Again this formula makes sense provided that we do not have $a=1$ and $a'=1$ simultaneously. But, in that case both elements are the unit of $A\t B$ for which we do not need to define the left product with itself.
  
  It is straightforward to show that, equipped with these two products, $A\t B$ becomes a unital dipterous algebra (cf. \cite{LRnote}).

 By definition a \emph{dipterous bialgebra} is a unital dipterous algebra $(\HH, *, \d)$ equipped with a coassociative counital coproduct $\DD: \HH \to \HH\t \HH$ which is a morphism of unital dipterous algebras. 
  This last condition means that  the following two compatibility relations hold:
  $$\DD(x *y)= \DD(x)*\DD(y) \quad \textrm{and} \quad \DD(x\d y)= \DD(x)\d \DD(y)= \sum x_{(1)}*y_{(1)}\t x_{(2)}\d y_{(2)}.$$

\begin{prop}\label{cofreedipterous} Any cofree-coassociative CHA  has a natural structure of unital dipterous bialgebra.
 \end{prop}
 \begin{proo} We denote by $*$ the associative product on $T^c(R)$ induced by the ${\mbox {\it MB}}$-structure of $R$. Define a new binary operation on 
 $T^c(R)$  by:
 $$(u_1\dots u_k)\succ (v_1\dots v_l):=((u_1\dots u_k)*(v_1\dots v_{l-1}))\succ v_l,$$
 for $u_1,\dots ,u_k,v_1,\dots ,v_l\in R$ for $l\geq 2$ and for $l=1$:
  $$(u_1\dots u_k)\succ v_1=u_1\dots u_k v_1.$$
In particular, we get:
 $$(((u_1\succ u_2)\succ u_3)\succ \dots )\succ u_k:=u_1u_2u_{3}\dots u_k.$$
 
 It is immediate to verify that the two binary operations $*$ and $\d$ satisfy the dipterous relation.
  
  The compatibility relation with the deconcatenation $\Delta $ is also immediate by induction.
    \end{proo}

    \begin{prop}\label{dipterousbialgebra} \cite{LRnote} The free unital dipterous algebra on the vector space $V$, denoted $Dipt(V)$, is naturally equipped with a dipterous bialgebra structure, and, a fortiori, with a Hopf algebra structure.
  \end{prop}  
  \begin{proo}  We follow the general procedure given in \cite{JLLsci}. Consider the map $V\to Dipt(V)\t Dipt(V), v \mapsto v\t 1 + 1 \t v$. Since $Dipt (V)$ is free, there exists a unique extension of this map as a dipterous morphism
  $$\DD: Dipt(V) \to  Dipt(V)\t Dipt(V).$$
  Using again the universal property of $Dipt(V)$ it is immediate to check that $\DD$ is coassociative. Hence $(Dipt(V), *, \d , \DD)$ is a dipterous bialgebra and this structure is functorial in $V$.
  \end{proo}
  
\subsection{Remark} Let $\HH$ be an $As^c-Dipt$-bialgebra which is conilpotent. By results of \cite{LRNote} there exists an isomorphism of coassociative coalgebras $T^c(\Prim \HH ) \cong \HH$. Any choice of such an isomorphism makes $\HH$ into a cofree-coassociative CHA.

  \subsection{From dipterous algebras to multibrace algebras}\label{dipttoMB} Let $A= (A, \succ, *)$ be a dipterous algebra. For elements $u_1,\dots ,u_p$ in $A$ we define :
  $$\omega ^{\succ }(u_1\dots u_p):= (((u_1\succ u_2)\succ u_3)\succ \dots )\succ u_p.$$
We construct $(p+q)$-ary operations on $A$ inductively as follows:
$$M_{00}= 0, \ M_{10}= \id_A = M_{01}, \ {\rm and }\ M_{n0}=0=M_{0n}\ {\rm
for} \ n\geq 2,$$
and
$$
M_{pq}:=\
 \big(\omega^{\d}\big) * \big(\omega^{\d}\big)
\hfill -\sum_{k\geq 2}\sum \omega^{\d}\big(M_{i_1j_1}M_{i_2j_2}\cdots
M_{i_kj_k}\big)$$
 where
the second sum (for which
$k\geq 2$ is fixed) is extended to all sets of indices $(i_1, \cdots,
i_k ; j_1,
\cdots , j_k)$ such that $i_1+ \cdots+ i_k = p$ and  $j_1+ \cdots+j_k = q$.

For instance:
$$\displaylines{M_{11}(u; v)= u*v -  u\d v - v\d u\ ,\hfill\cr
}$$
$$\displaylines{
M_{21}(uv ; w)= (u\d v)*w   - u\d M_{11}(v; w) -   M_{11}(u; w)\d v\hfill\cr
\hfill  -(u\d v)\d w  - (u\d w)\d v  - (w\d u)\d v\ ,  \cr
\qquad = (u\d v)*w -(u\d v)\d w - u\d (v*w) + u\d(v\d w),\cr
}$$
$$\displaylines{
M_{12}(u; vw)= u*(v\d w)  - M_{11}(u; v)\d w -  v\d M_{11}(u; w)\hfill\cr
 \hfill  -(u\d v)\d w  -  (v\d u)\d w  -  (v\d w)\d u\ . \cr
 \qquad= u*(v\d w) -(u*v)\d w - v\d (u*w - u\d w - w\d u) - (v\d w)\d u.\cr
}$$

Let us remark that, if we define $x\g y:=x*y - x\d y$, then $M_{21}$ becomes:
$$M_{21}(uv ; w)=  (u\d v)\g w -  u\d (v\g w).$$

 \begin{prop}\label{functorFcase1} There is a (forgetful) functor $F$ from the category of dipterous algebras to the category of multibrace algebras given by 
  $$(A, *, \d) \mapsto (A, \{M_{pq}\}_{p\geq 0, q\geq 0}).$$
  The functor $F$ has a left adjoint denoted $U$. For a brace algebra $R$, the dipterous algebra $U(R)$  is the quotient of the free dipterous algebra on the vector space $R$ modulo the relations which identify the $MB$-structure of $R$ with the $MB$-structure coming from the dipterous structure. Moreover $U(R)$ is a dipterous bialgebra.
  \end{prop}

\begin{proo}  The proof follows immediately from the definition of a multibrace algebra, using that 
$\omega ^{\succ}(u_1\dots u_p)=u_1\otimes \dots \otimes u_p$ and that $*$ is associative.

Since $F$ is a forgetful functor, it is immediate that it has a left adjoint, which is precisely $U$ as described in the statement.

Since the free dipterous algebra is a dipterous bialgebra by Proposition \ref{dipterousbialgebra}, it follows that its quotient $U(R)$ is also a dipterous bialgebra.
\end{proo}

\medskip

\begin{thm}\label{tripleAsDiptMB} There is a good triple of operads  
$$(As, Dipt, MB).$$
For any vector space $V$ there is a natural isomorphism of CHAs:
$$Dipt(V) \cong T^c(MB(V)).$$
\end{thm}
\begin{proo}
Let $\PP= \Prim_{As}Dipt$ denote the operad made of the primitive operations for $As^c\pt Dipt$-bialgebras. The hypothesis (H1), recalled in \ref{triples}, is fullfilled thanks to Proposition \ref{dipterousbialgebra}. The hypothesis (H2epi) is also fullfilled because a splitting of the coalgebra map $Dipt(V)\to As(V)=T^c(V)$ is given by
$$v_{1}\cdots v_{n}\mapsto \omega^{\d}(v_{1}\cdots v_{n}).$$
The coalgebra morphism property is proved by induction on $n$.  Hence, by \cite{JLLgbto}, there is a good triple of operads
$$(As, Dipt, \PP).$$
It remains to show that we have $\PP= MB$.  Since the operations $M_{pq}$ are primitive by construction there is a natural map $MB(V) \to \PP(V)$, with extension
$T^c(MB(V)) \to T^c(\PP(V))$. By Theorem \ref{coAsCombHopfMB} $T^c(MB(V)) $ is a dipterous algebra, therefore there is a natural dipterous morphism from $Dipt(V)$ to it. Since $(As, Dipt, \PP)$ is a good triple of operads, the composite 
$$Dipt(V) \to T^c(MB(V)) \to T^c(\PP(V))$$
is an isomorphism. As a consequence  the map $MB(n)\to \PP(n)$ on the spaces of $n$-ary operations is surjective. It suffices to show that they have the same (finite) dimension to deduce that it is an isomorphism.

The dimension of $MB(n)$ is known to be $C_{n} \times n!$, where $C_{n}$ is the 
super-Catalan number (number of planar rooted trees with $n$ leaves), cf.\ \cite{LRstr}. The dimension of $Dipt(n)$ is known to be $2\times C_{n} \times n!$, see \ref{freediptalg}. From the functional equation relating the generating series of the two functors $Dipt$ and $\PP$ deduced from $Dipt = As \circ \PP$  and from the functional equation satisfied by the generating series of the super-Catalan numbers (cf.\ \ref{planartrees}) we conclude that $\dim \PP (n)=C_{n} \times n!$. Hence $\dim MB(n) = \dim \PP (n)$ and $MB=\PP$ as expected. 
\end{proo}

\subsection{Planar trees}\label{planartrees}

In order to describe the structure of the free dipterous algebra we introduce the combinatorial objects named 
planar trees.  We denote by $PT_n$ the set of \emph{planar
trees} with $n$ leaves, $n\geq 1$, which are reduced (every vertex has at least two inputs) and rooted.
Here are the first of them:
$$PT_1 = \{\vert\} ,\qquad \ PT_2 = \{\ \arbreA\} ,\qquad PT_3 = \Big\{\
\arbreAB , \arbreBA , \arbreBB\Big\}.$$

The lowest vertex is called the {\it root vertex}. Observe that the tree $\mid$ has no vertex.

We define $PT_{\infty} := \bigcup_{n\geq 1} PT_n$. The number of elements
in $PT_n$ is the so-called {\it super
Catalan number} $C_n$. 

By definition the {\it grafting} of $k$ planar trees $\{x^{(1)}, \ldots ,
x^{(k)}\}$ is a planar tree denoted $\bigvee (x^{(1)},
\cdots, x^{(k)})$ obtained by joining the $k$ roots to a new vertex and adding
a new root. For $k=2$, sometimes we shall denote $x^{(1)}\vee x^{(2)}$ instead of $\bigvee (x^{(1)},x^{(2)})$. 
Observe that the grafting operation is not associative. In fact the three trees $(x\vee y)\vee z,
x\vee (y\vee z), \bigvee(x, y, z)$ are all
different.
Any planar tree $x$ can be uniquely obtained as
$x=\bigvee (x^{(1)},\cdots ,x^{(k)})$, where $k$ is the number of inputs of the root
vertex. 

The grafting operation gives a bijection between the set of $k$-tuples of
planar trees, for any $k\geq 2$, and
$PT_{\infty}\backslash \{\mid\}$. Indeed the inverse is the map
$\bigvee (x^{(1)},\cdots ,x^{(k)})\mapsto\{x^{(1)}, \cdots , x^{(k)}\}$. As a consequence the generating series $C(t) = \sum_{n\geq 0} C_{n}t^n$ satisfies the equation $2tCf(t)^2 - (1+t)C(t) +1 = 0$.

\subsection{Free dipterous algebra}\label{freediptalg} Since, in the relations defining the
notion of a dipterous algebra, the variables stay in the same
order, we only need to understand the free dipterous algebra on one generator (i.e.\ over $\KK$). In other words the operad $Dipt$ is a nonsymmetric operad (see \cite{LV}).

Let $T(\PTi)$ be the free unital associative algebra over the vector space
$\KK[\PTi]$ spanned by the set $\PTi$. So
we write $T(\PTi)$ instead of $T(\KK[\PTi])$. A set of linear generators of
$T(\PTi)$ is made of the monomials
$t_1 t_2 \cdots t_k$ where the $t_i$'s are planar trees. We
define a right product on the augmentation ideal as follows:
$$(t_1 t_2 \cdots t_k)\d (s_1 s_2 \cdots s_l) := \Big( t_1 \vee \big(t_2
\vee( \cdots \vee \bigvee (t_k,s_1, \ldots , s_l)\cdot )\big)\Big).$$
Observe that the right product of two monomials is always a tree. As before
we extend the right operation  by $1\succ
\omega= \omega$ and $\omega\succ 1=0$ for $\omega$ in the augmentation ideal.
\medskip

\begin{prop}\label{freediptandtrees} The associative algebra $T(\PTi)$ equipped
with the right product defined above is the free 
unital dipterous algebra on one generator.\end{prop}
\medskip

\begin{proo} Let us verify the dipterous relation. On one hand we have
$$\displaylines{
((t_1   \cdots t_k)*(s_1   \cdots s_l))\d (r_1   \cdots r_m) = (t_1
\cdots t_k s_1   \cdots s_l)\d (r_1
\cdots r_m) \hfill \cr
\hfill = \Big(t_1\vee \big(   \cdots \vee(t_k \bigvee (s_1,   \cdots , s_{l-1},
\bigvee(s_l, r_1,\cdots , r_m) \cdot ))\cdot \big)\Big)\cr
}$$
On the other hand we have
$$\displaylines{
(t_1   \cdots t_k)\d \big((s_1   \cdots s_l)\d (r_1   \cdots r_m)\big) =
(t_1   \cdots t_k)\d \big( s_1\vee(   \cdots
\vee (s_l\vee r_1
\cdots \vee r_m)\cdot \big)\Big) \hfill
\cr
\hfill= \Big(t_1\vee \big(   \cdots \vee(t_k \bigvee (s_1,\cdots ,s_{l-1},\bigvee (s_l, r_1
\cdots , r_m) \cdot ))\cdot \big)\Big)\ .\cr
}$$
And so $(T(PT_{\infty}, *, \d)$ is a dipterous algebra. Since any linear generator can be obtained from the tree $|$ by combining the two operations $*$ and $\d$, it is a quotient of $Dipt(\KK)$. In order to show that it is free on one generator it is sufficient to show that, for any dipterous algebra $A$ and any element $a\in A$, there is a unique dipterous morphism $T(PT_{\infty})\to A$ mapping $|$ to $a$. We already know how to construct a map and we know that it is unique. To prove that it is a dipterous morphism we proceed by induction (on the degree of the trees) like in the proof of the dendriform case performed in \cite{JLLdig} Proposition 5.7.\hfill \end{proo}

\subsection{Remark} As in the case of the operad $2as$ governing the 2-associative algebras (cf.\ \cite{LRstr}) the operad $Dipt$ can be described in terms of two copies of the set $PT_{\infty}$. 
  
  \subsection{Free-associative CHA and co-multibraces}\label{freeassCHA}
  
A \emph{free-associative combinatorial Hopf algebra} $\HH= \KK 1\oplus \HHo$, free-associative CHA for short, is a Hopf algebra structure on the tensor algebra $T(C)$. Here $C$ is the space of indecomposables of $\HH$, that is $C:= \HHo/ \HHo^2$. This structure is encoded into the coproduct map $\DD$. The preceding results admit an obvious dualization. So, $C$ is a \emph{multibrace coalgebra} and this comultibrace structure determines completely the Hopf algebra $(T(C), *,\DD)$. The space $T(C)$ inherits a \emph{dipterous coalgebra} structure which makes it into a $Dipt^c \pt As$-bialgebra.

For any vector space $V$ there is an isomorphism of CHAs 
$$T(MB^c(V))\cong Dipt^c(V).$$

These spaces are spanned by trees and the (co)operations can be made explicit through ``admissible cuts'' (the dual notion of grafting).

  \subsection{2-associative bialgebras}\label{2as} A 2-associative algebra $A$ (or $2as$-algebra for short) is a vector space equipped with two associative products. Here we suppose that they are unital, that is $A$ contains an element $1$ which is a unit for the two associative products $\cdot$ and $*$. In \cite{LRstr} we introduced the notion of 2-associative bialgebra (more precisely $As^c\pt 2as$-bialgebra) which is a $2as$-algebra equipped with a coassociative (and counital) coproduct satisfying the Hopf compatibility relation with $*$ and the unital infinitesimal compatibility relation with $\cdot$ (see loc.cit.\ for details).
  
  Let $\HH=(T^c(R), *)$ be a cofree-coassociative CHA. In \cite{LRstr} we have shown that, if we consider the concatenation product $\cdot$ on the tensor module $T^c(R)$, then $(T^c(R), *, \cdot, \DD)$ is a 2-associative bialgebra. The same results as before are valid with dipterous replaced by 2-associative. They are the subject of \cite{LRstr} where it is proved that 
$$(As, 2as, MB)$$
 is a good triple of operads. 
 
 Let us recall from loc.\ cit.\ that there is a functor
 $$F: 2as\alg \to MB\alg$$
 which provides operations $M_{pq}$ out of the two associative operations $\cdot$ and $*$.
For instance:
$$\displaylines{
M_{11}(x; y)= x*y -  x\cdot y - y\cdot x\ ,\hfill\cr
M_{21}(xy ; z)= (x\cdot y)*z   - x\cdot (y*z) - (x*z)\cdot y + x\cdot y \cdot z ,\hfill \cr
M_{12}(x; yz)= x*(y\cdot z) -(x*y)\cdot z - y\cdot (x*z)+  y\cdot x \cdot z.\hfill\cr
}$$

 The choice of 2-associative bialgebras behaves well with the symmetrization procedure that we deal with in Section \ref{s:cofreecocom}. The choice of dipterous bialgebras behaves well with the right-sided hypothesis, dealt with in Section \ref{s:cofreecoasrs}.
  
\subsection{From dipterous to 2-associative}\label{Dipt2as} Let $\lambda$ be a formal parameter (or an element in $\KK$). We consider algebras having two binary operations $*$ and $\d$ such that $*$ is associative and the following relation holds:
$$(\lambda x*y + (1- \lambda) x\d y )\d z = x \d (y\d z).$$
Let us denote by $\ZZZ_{\lambda}$ the associated operad. It is clear that $\ZZZ_{0}= 2as$ and $\ZZZ_{1}=Dipt$. So we have a homotopy between the operads $2as$ and $Dipt$. It would be interesting to know if we always have a good triple of operads
$$(As, \ZZZ_{\lambda}, MB)$$
such that any cofree-coassociative CHA is a $\ZZZ_{\lambda}$-bialgebra.

 \section{(Co)free-(co)associative right-sided CHA}\label{s:cofreecoasrs} 
 
 In this section we study the cofree-coassociative CHAs (and then free-associative CHAs) which satisfy the right-sided property (r-s). We will see that they are dendriform bialgebras which have been thoroughly studied in \cite{JLLdig, LRHopf, Ronco00, Ronco01, Ronco02}. The space of primitives inherits a brace algebra structure. In the last part of the section we study the alternative structure, denoted $\YYY$, which is close to the 2-associative structure. Examples will be given in section \ref{Exs}. They include  the Hopf algebra of quasi-symmetric functions (\ref{qsym}), Malvenuto-Reutenauer algebra (\ref{MalvenutoReutenauer}), the Solomon-Tits algebra (\ref{qsym}), the algebra of parking functions (\ref{qsym}) and  Hopf algebras constructed out of operads (\ref{operads}).
 
 \subsection{Right-sided condition on combinatorial Hopf algebras}\label{conditionr-s} Let $\HH=(T^c(R), *)$ be a cofree-coassociative CHA. There is a  natural grading  $\HH=\oplus_{n} \HH^n$ given by $\HH^n:= R^{\t n}$. We study the cofree-coassociative CHAs satisfying the following  condition.
 
 \medskip
 
 \noindent Right-sided condition:
 
 \noindent {\bf (r-s)} for any integer $q$ the subspace $\bigoplus_{n\geq q} {\HH^n}$ is a right-sided ideal of $\HH$. 
 
 \medskip
Since $\HH$ is a cofree-coassociative CHA, by Proposition \ref{Rijk} the associative product $*$ is given by a family 
$\{M_{pq}\}_{p,q\geq 0}$ of  $(p+q)$-ary operations, satisfying the relations $\RRR$.

\begin{prop}\label{rsbrace} The combinatorial Hopf algebra $\HH$ is  right-sided if, and only if, $M_{pq}=0$ for $p\geq 2$.\end{prop}

\begin{proo} Recall that $\HH= T^c(V)$ is right-sided if the subspace ${\mathbb F}_r=\bigoplus _{n\geq r}V^{\t n}$ 
 is a right ideal of $T^c(V)$. 

If ${\mathbb F}_r$ is a right-sided ideal, then 
$$(x_1\dots  x_p) *(y_1\dots  y_q)=M_{pq}(x_1\dots x_p;y_1\dots y_q) + z,$$
 where $z\in {\mathbb F}_2$. But, for $p\geq 2$, $(x_1\dots  x_p) *(y_1\dots  y_q)$ must belong to ${\mathbb F}_p$ and $M_{pq}(x_1\dots x_p;y_1\dots y_q)\in {\mathbb F}_1\setminus {\mathbb F}_2$, which implies that $M_{pq}=0$ for all $q$ and $p\geq 2$.
\medskip

Conversely, suppose that $M_{pq}=0$  for all $m$ and $n\geq 2$. 
The formula for $*$ is then given by:
$$\begin{array}{l}
(x_1\dots  x_p) *(y_1\dots  y_q)=\\
{\displaystyle \sum _{r}\big(\sum _{{{\underline k}\vdash q}\atop {|k|=p+2}}}y_1\dots y_{k_0}M_{1k_1}(x_1;y_{k_0+1},\dots ,y_{k_0+j_{1}})\dots 
y_{k_{1}}M_{1j_2}(x_2;\dots )\dots y_{k_{1}+j_{2}}) \dots \\
\hfill M_{1j_p}(x_p;\dots )\dots y_{q})\big),
\end{array}
$$
where $ 0\leq k_0\leq k_{0}+j_{1}\leq \dots \leq  k_{p-1}+j_{p}\leq  q$.
So, the element $(x_1\dots  x_p) *(y_1\dots  y_q)$ belongs to ${\mathbb F}_p$, which ends the proof.
\end{proo}

 \subsection{Brace algebras}\label{brace} A \emph{brace algebra} is a vector space $R$ equipped with  a $(1+q)$-ary operation denoted $\{-; -, \ldots , -\}$ for any $q\geq 1$ satisfying the following formulas:
 $$ \{\{x;y_1,\ldots , y_n\};z_1,\ldots ,z_m \} = \sum\{x; \ldots ,\{y_1; \ldots\}, \ldots , \{y_n; \ldots , \}, \ldots \}.$$
 On the right-hand side the dots are filled with the variables $z_i$'s (in order) with the convention $\{y_k; \emptyset \}= y_k$. This notion appears in the work of 
 Gerstenhaber and Voronov in \cite{GerstenhaberVoronov}.
 
 \begin{lemma}\label{multibracetrivial} A multibrace algebra for which $M_{pq}=0$ when $p\geq 2$ is a brace algebra and vice-versa.
 \end{lemma}
 \begin{proo} It is straightforward to check that, under the identification
 $$M_{1q}(x; y_{1}\cdots y_{q})= \{x; y_{1},\ldots , y_{q}\},$$
 the relations $\RRR_{ijk}$ are exactly the brace relations.
 \end{proo}
 
 \begin{thm} The primitive part of a cofree-coassociative r-s CHA  is a brace algebra.  There is an equivalence of categories between cofree-coassociative r-s CHAs and brace algebras.
 \end{thm}
 \begin{proo} By Proposition \ref{Rijk} the primitive part of a cofree bialgebra is a $MB$-algebra. By Lemma \ref{rsbrace} property (r-s) implies that $M_{pq}=0$ when $p\geq 2$. Therefore, by Lemma \ref{multibracetrivial} this primitive part is a brace algebra. The two functors $R\mapsto T^c(R)$ and $\HH = T^c(R) \mapsto R$ are obviously  inverse to each other.
 \end{proo}
 
 \subsection{Dendriform algebras \cite{JLLdig}}\label{dend} A \emph{dendriform algebra} is a vector space $A$ equipped with two binary operations denoted $\g$ and $\d$ which satisfy the conditions:
 \begin{displaymath}
\left\{
\begin{array}{ccc}
(x\g y)\g z &=& x\g (y* z)  , \\
(x\d y)\g z &=& x\d (y\g z) , \\
(x* y)\d z &=& x\d (y\d z).
\end{array}
\right.
\end{displaymath}
 where $x*y:= x\g y + x \d y$.
 
 As in the case of dipterous algebras, there is a notion of \emph{unital dendriform algebra}. It is a vector space $A=\KK \, 1 \oplus \Ao$ such that $\Ao$ is a dendriform algebra and the two binary operations are extended as follows:
 \begin{displaymath}
\left\{\begin{array}{rclcrcl}
1\g x &=& 0 & , & x\g 1& = &x ,\\
1\d x &=& x & , & x\d 1& =& 0 .\\
\end{array}
\right.
\end{displaymath}
Observe that one has $1*x=x=x*1$ as expected, but $1\g 1$ and $1\d 1$ are \emph{not} defined.

\begin{prop}\label{CHAtodend} Any cofree-coassociative CHA, which is right-sided, is a dendriform bialgebra.
 \end{prop}
 \begin{proo} By Proposition \ref{cofreedipterous} the Hopf algebra $\HH$ is a dipterous algebra. Let us define the operation $\g$ by the equality $x*y= x\g y + x \d y$.
 
 Condition (r-s) implies that the operations $\g$ and $\d$ satisfy the relation 
 $$(x\d y)\g z = x\d (y\g z).$$
 It is clear that the axioms for a dipterous algebra satisfying this extra condition are equivalent to the axioms for dendriform algebras.
 \end{proo}
 
 \subsection{From dendriform algebras to brace algebras}\label{dendto brace} The forgetful functor from dipterous algebras to multibrace algebras sends a dendriform algebra to a brace algebra with the same underlying space. In  fact its restriction is precisely the functor constructed in \cite{Ronco00}. Explicitly it is given by
 $$M_{1q}(v; v_{1}, \ldots ,v_{q})= 
 \sum_{i=0}^{i=q}(-1)^{i}\omega_{\g}(v_{1} \ldots v_{i})\d v \g \omega^\d(v_{i+1} \ldots v_{q}),$$
where $$\omega_{\g}(v_{1} \ldots v_{n}):= v_1\g(v_2\g(\dots \g (v_{n-1}\g v_n)))$$ and where $$ \omega^\d(v_{1} \ldots v_{n}):= ((v_1\d v_2)\d \dots )\d v_n,$$
is the element defined in Subsection 2.9. 
 The proofs of the following results are similar to those of section \ref{s:cofreecoas} and can be found in \cite{Ronco00, Ronco01, Ronco02}.
 
 \begin{thm} There is a good triple of operads
 $$(As, Dend, Brace).$$
  For any vector space $V$ there is a natural isomorphism of cofree dendriform bialgebras, and, a fortiori, of cofree-coassociative right-sided CHAs:
  $$Dend(V)\cong T^c(Brace(V)) .$$
 \end{thm}

\subsection{Planar binary trees and free dendriform algebra}\label{pbtDend} The set of planar binary trees with $n$ leaves, denoted $PBT_{n}$, is made of the planar trees (cf.\ \ref{planartrees}) with exactly two inputs at each vertex. The number of elements of $PBT_{n}$ is the Catalan number $c_{n}$. For any $n\geq 1$ a tree $t\in PBT_{n}$ can be uniquely written as a grafting $t= t^l\vee t^r$. It is shown in \cite{JLLdig} that the free dendriform algebra on one generator $Dend(\KK)$ is the Hopf algebra $\HH_{LR}:=\oplus_{n\geq 1}\KK[PBT_{n}]$ equipped with the following right and left product:
$$t\g s := t^l \vee (t^{r}*s), \quad t\d s := (t*s^l)\vee s^r.$$
It was first introduced in our paper \cite{LRHopf}. 
Notice that the tree $|$ is the unit and the tree $y:=\arbreA$ is the generator. In order to compute in $\HH_{LR}$ the following formula proves to be helpful:
$$t\d y \g s = t \vee s.$$
In order to work with the free dendriform algebra on a vector space $V$ (resp. on a set $X$) it suffices to decorate the trees by putting  elements of $V$ (resp. $X$) in between the leaves. For instance:
$$\arbreAv , \arbreABv, \arbreBAv, \arbreACAv,$$
represent respectively
$$v ,\quad  v_{1}\d v_{2}, \quad  v_{1}\g v_{2}, \quad  v_{1}\d v_{2}\g v_{3}.$$
 
  \subsection{Planar binary trees and the primitives}\label{pbtprimitives} The primitive operations can be described by planar binary trees as follows. The following results, stated without proofs, are taken from \cite{Ronco02}.
 
The dendriform analog of the Eulerian idempotent acting on the dendriform bialgebra $Dend(\KK)$ is defined by the formula 
$$e:=\sum _{n\geq 1}(-1)^{n+1}\omega ^n_{\d}\circ {\overline {\Delta }}^n,$$
where ${\overline {\Delta }}$ is the reduced coproduct of  $Dend(\KK)$, and 
where 
$$\omega _{\succ }^n(x_1,\dots ,x_n):=x_1\succ (x_2\succ (\dots \succ (x_{n-1}\succ x_n))), \textrm{ for } x_{i}\in Dend(\KK).$$
This operator has the following properties:

\begin{enumerate}
\item the element $e(x)$ is a primitive element, for any $x\in  Dend(\KK)$,
\item $e(x\succ y)=0$ , for any pair of elements $x$ and $y$ of $Dend(\KK)$,
\item $e(x)=x$ for any $x\in {\rm Prim}(Dend(\KK))$.
\end{enumerate}

Moreover, if $t=\vert \vee x$ is a tree in $PBT_n$, with $x\in PBT_{n-1}$, then $e(t)=t+\sum_{i} t_1^i\vee t_2^i$, for some  $t_1^i\in PBT_{n_i}$ with $n_i\geq 2$.

An easy argument on the dimensions of $Dend(\KK)_n$ and $T(\bigoplus_{m\geq 0}\KK[ PBT_m])_n$ shows that the elements $e(\vert \vee t)$, with $t\in PBT_n$, form a basis of 
${\rm Prim}(Dend(\KK))_n$, for all $n\geq 1$. 
From these results it follows that any element of $Dend(\KK)$ is a sum of elements of type 
$((x_1\succ x_2)\succ \dots )\succ x_r$, with $x_i\in {\rm Prim}(Dend(\KK))$ for $1\leq i\leq r$. 

Moreover, the coproduct $\Delta $ of $Dend(\KK)$, satisfies that:
$$\Delta (((x_1\succ x_2)\succ \dots )\succ x_r)=\sum _{i=0}^r((x_1\succ x_2)\succ \dots 
)\succ x_i\otimes ((x_{i+1}\succ x_{i+2})\succ \dots )\succ x_r,$$
for any $x_1,\dots ,x_r\in {\rm Prim}(Dend(\KK))$.
Therefore there is a bijection 
$$\KK[PBT_{n-1}\times S_{n}]\cong (\Prim_{As}Dend)(n)$$
 given by
$(t;x_{1},\ldots , x_{n}) \mapsto (e(\mid\vee t); x_{1},\ldots ,
x_{n})$. For instance, in the free dendriform algebra $Dend(V)$ over the vector space $V$ we get the following formulas in low dimension for any $v_{i}\in V$:
$$\displaylines{
(\mid; v_{1}) \mapsto v_{1},\quad  (\arbreA; v_{1} v_{2})\mapsto v_{1}\g
v_{2}- v_{2}\d v_{1}= M_{11}(v_{1}; v_{2}),}
$$
$$\displaylines{
\Big(\arbreAB; v_{1}v_{2}v_{3}\Big)\mapsto\hfill \cr
\qquad \qquad \Big(\arbreCAB; v_{1}v_{2}v_{3}\Big)-\Big (\arbreACA; 
v_{2}v_{1}v_{3}\Big)+\Big (\arbreBAC; v_{2}v_{3}v_{1}\Big)\hfill\cr
\qquad \qquad =v_{1}\g (v_{2}\d v_{3})- v_{2}\d v_{1}\g v_{3})+ (v_{2}\g
v_{3})\d v_{1})\hfill\cr
\qquad \qquad = M_{12}(v_{1};v_{2}v_{3}) ,\hfill\cr
}$$
$$\displaylines{
\Big(\arbreBA; v_{1}v_{2}v_{3}\Big)\mapsto\hfill \cr
\qquad \Big(\arbreCBA; v_{1}v_{2}v_{3}\Big)-\Big(\arbreACA; v_{3}v_{1}v_{2}\Big)-\Big(\arbreBAC; 
v_{2}v_{3}v_{1}\Big)\hfill \cr
\hfill +\Big(\arbreABC;v_{3}v_{2}v_{1}\Big)+\Big(\arbreBAC; 
v_{3}v_{2}v_{1}\Big)\cr
\qquad =v_{1}\g (v_{2}\g v_{3})- v_{3}\d v_{1}\g v_{2}- (v_{2}\g v_{3})\d v_1+ (v_{3}\d
v_{2})\d v_{1}+ (v_{3}\g v_{2})\d v_{1}\hfill \cr
\qquad = M_{11}(M_{11}(v_{1};v_{2});v_{3})- M_{12}(v_{1};v_{2}v_{3}).\hfill \cr
}$$

Various basis of $Dend(\KK)$ and comparison between them have been studied in \cite{AguiarSottile, Foissy, Holtkamp03}.

\subsection{Unreduced planar trees}\label{uptrees} We denote the set of \emph{planar 
unreduced trees} with $n$ vertices by $PUT_{n}$. In low dimensions we have, for $n=1$ to 
$4$,

{\unreducedtrees}

The number of elements of $PUT_{n}$ is also the Catalan number $c_{n}$. Let us define a grafting operation on $PUT_{\infty}$ as follows. The unreduced planar tree $x\vee y$ is obtained by putting the tree $y$ on the right-hand side of $x$ and joining the root of $x$ to the root of $y$ by a new edge. The root-vertex of $x\vee y$ is taken to be the root-vertex of $x$. 

\begin{lemma}\label{PUTvsPBT}\emph{ [P.\ Palacios, thesis, unpublished].} The map $\varphi : PBT_{n}\to PUT_{n}$ given by $\varphi(t\vee s) = \varphi(t) \vee \varphi(s)$ and $\varphi(|) = \bullet$ is a bijection. 
\end{lemma}

The following list of planar binary trees gives under $\varphi$ the list of planar unreduced trees cited above:
$$\mid \arbreA \arbreBA \arbreAB \arbreCBA \arbreCAB \arbreACA \arbreBAC \arbreABC$$

\subsection{Labelled planar unreduced trees and the operad $Brace$}\label{lput} 
We recall the explicit description of the operad $Brace$ in terms of planar unreduced trees (see for instance \cite{Chapoton02}). We first construct the $S_{n}$-module $Brace(n)$ and then we give the formula for the partial composition operation $\circ_{i}$.

The $S_{n}$-module $Brace(n)$ is the following regular representation:
$$Brace(n)= \KK[PUT_{n}\times S_{n}.$$
The generating operation $M_{1n}$ corresponds to $c_{n}\times \id_{n}$ where $c_{n}$ is the corolla. In order to describe $\mu\circ_{i}\nu$ for $\mu\in Brace(m) $ and $\nu\in Brace(n)$ it suffices to describe it for $\mu= t\times \id_{m}$ and $\nu=\r\times \id_{n}$ ($t$ and $r$ are trees). The composite
$$( t\times \id_{m}) \circ_{1}( r\times \id_{n})$$
is obtained through the brace relation \ref{brace}. For instance, letting $n=m=2$, the relation
$$\{\{x_{1}; x_{2}\};x_{3}\} = \{\{x_{1}; x_{2},x_{3}\} + \{x_{1}; \{x_{2};x_{3}\}\} +  \{x_{1}; x_{3};x_{2}\} $$
gives:

\begin{picture}(400,50)
\put(30,15){\scriptsize $\bigl ($}
\put(35,5){\scriptsize$\bullet$}
\put(35,27){\scriptsize $\bullet$}
\put(37,12){\line(0,1){12}}
\put(40,15){\scriptsize $\times$}
\put(47,15){\scriptsize $\id_2$}
\put(62,15){\scriptsize $\bigr )$}
\put(67,15){\scriptsize $\circ _1$}
\put(74,15){\scriptsize $\bigl ($}
\put(79,5){\scriptsize$\bullet$}
\put(79,27){\scriptsize $\bullet$}
\put(81,12){\line(0,1){12}}
\put(84,15){\scriptsize $\times$}
\put(92,15){\scriptsize $\id_2$}
\put(107,15){\scriptsize $\bigr )$}
\put(112,15){\scriptsize $=$}
\put(119,24){\scriptsize$\bullet$}
\put(127,5){\scriptsize$\bullet$}
\put(135,24){\scriptsize$\bullet$}
\put(127,7){\line(-1,3){6}}
\put(131,7){\line(1,3){6}}
\put(140,15){\scriptsize $\times $}
\put(150,15){\scriptsize $\id_3 $}
\put(163,15){\scriptsize $+$}
\put(170,5){\scriptsize$\bullet$}
\put(170,27){\scriptsize $\bullet$}
\put(170,49){\scriptsize $\bullet$}
\put(172,12){\line(0,1){12}}
\put(172,34){\line(0,1){12}}
\put(177,15){\scriptsize $\times $}
\put(187,15){\scriptsize $\id_3 $}
\put(200,15){\scriptsize $+$}
\put(207,24){\scriptsize$\bullet$}
\put(215,5){\scriptsize$\bullet$}
\put(223,24){\scriptsize$\bullet$}
\put(215,7){\line(-1,3){6}}
\put(219,7){\line(1,3){6}}
\put(228,15){\scriptsize $\times $}
\put(238,15){\scriptsize $[132]$}
\end{picture}

Observe that the action of the symmetric group is involved in this formula. In order to describe the $\circ_{i}$ composition operation we label the vertices of a planar unreduced tree as follows. The root-vertex is labelled by $1$. The other vertices are labelled according to the following rule: for $t=x\vee y$ the integers labelling the vertices of $x$ are less than the the integers labelling the vertices of $y$.  Here are the first of them:

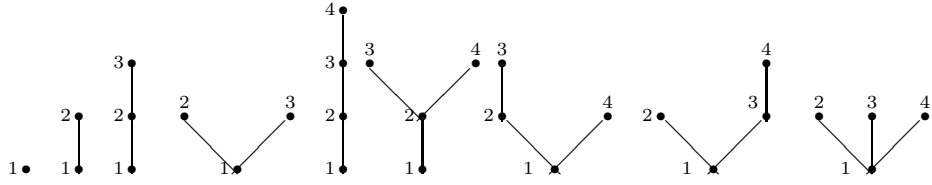
\labeledunreducedtrees

For $i=2, \ldots , m$ we have 
$$( t\times \id_{m}) \circ_{i}( r\times \id_{n})= (t\circ_{i}r \times \id_{m+n-1}$$
where the $t\circ_{i}r $ is obtained by grafting the tree $r$ to the tree $t$ by adjoining an edge from the $i$th vertex of $t$ to the root of $r$. For instance:

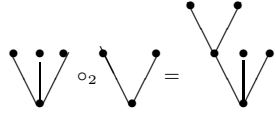
\begin{center}
\begin{picture}(400,50)
\put(30,5){\scriptsize$\bullet$}
\put(20,24){\scriptsize $\bullet$}
\put(30,24){\scriptsize $\bullet$}
\put(39,24){\scriptsize $\bullet$}
\put(31,6){\line(-1,2){10}}
\put(32,6){\line(0,1){16}}
\put(33,6){\line(1,2){10}}
\put(46,15){\scriptsize $\circ_2$}
\put(54,24){\scriptsize $\bullet$}
\put(65,5){\scriptsize$\bullet$}
\put(74,24){\scriptsize$\bullet$}
\put(66,6){\line(-1,2){11}}
\put(67,6){\line(1,2){10}}
\put(79,15){\scriptsize $=$}
\put(87,42){\scriptsize$\bullet$}
\put(105,42){\scriptsize$\bullet$}
\put(96,24){\scriptsize $\bullet$}
\put(107,5){\scriptsize $\bullet$}
\put(107,24){\scriptsize $\bullet$}
\put(117,24){\scriptsize $\bullet$}
\put(108,6){\line(-1,2){10}}
\put(109,6){\line(0,1){17}}
\put(109,6){\line(1,2){10}}
\put(98,25){\line(-1,2){10}}
\put(98,25){\line(1,2){10}}
\end{picture}
\end{center}

\subsection{Comparison of $\Prim_{As}Dend$ and $Brace$} Putting together the results of the preceding paragraphs we see that there is a bijection
$$\KK[PBT_{n}\times S_{n}]\xrightarrow{\varphi} \KK[PUT_{n}\times S_{n}]\cong Brace(n)= \Prim_{As}Dend (n)\cong \KK[PBT_{n}\times S_{n}],$$
where the last isomorphism is the inverse of $t\mapsto e(\mid \vee t)$, cf.\ \ref{pbtprimitives}.
This composite is not the identity. 

\subsection{Free-associative right-sided CHA}\label{freeassrsCHA}
In the dual framework, that is for an associative CHA $\HH=T(C)$, where $C$ is a multibrace coalgebra, the right-sided condition reads as follows:

\noindent {\bf (r-s)} the coproduct $\DD$ on $T(C)$, given by $\DD(v) =v\t 1 + 1 \t v + \sum v_{(1)} \t v_{(2)}$, is such that
$$v_{(2)}\in C \textrm{ for any } v\in C.$$

This condition implies that $C$ is in fact a brace coalgebra and that $\HH$ is a dendriform coalgebra. 

All the results of this section can be dualized: $C$ is a brace coalgebra and $T(C)$ is a $Dend^c\pt As$-bialgebra. For instance the free-associative CHA $Dend^c(\KK)$ can be seen as a noncommutative variation of the Connes-Kreimer algebra. The coproduct can be described on planar binary trees and hence on planar unreduced trees by means of admissible cuts. Some of these results have been worked out explicitly in \cite{Foissy}.

\subsection{$\YYY$-algebras}\label{ss:YYY} If, instead of looking at the dipterous structure of a cofree-coassociative {\bf r-s} CHA, we look its 2-associative structure (cf.\ \cite{LRstr} and \ref{2as}), then the relevant quotient structure is more complicated: it is a $\YYY$-algebra, where $\YYY$ is the quotient of the operad $2as$ by the operadic ideal generated by the operations $M_{pq}$ for $p\geq 2$:
$$\YYY := 2as/\{M_{pq}\ |\ p\geq 2\}.$$
Since the operations $M_{pq}$ are primitive for $2as$-bialgebras, the notion of $\YYY$-bialgebra is well-defined. Let $\HH$ be a cofree-coassociative CHA which is right-sided. By Theorem 4.2 of \cite{LRstr} it is a $2as$-bialgebra. By Proposition 3.2 the operations $M_{pq}$ are $0$ on $\HH$ for $p\geq 2$. Therefore $\HH$ is a $\YYY$-algebra, and even a $\YYY$-bialgebra.

\begin{thm}\label{tripleYYY} There is a good triple of operads
$$(As, \YYY, Brace).$$
\end{thm}
\begin{proo} Since the operations $M_{pq}$ are primitive, by Proposition 3.1.1 of \cite{JLLgbto} there is a good triple of operads
$$(As, \YYY,\Prim_{As}\YYY).$$
Let us consider the free $\YYY$-algebra over $V$, denoted $\YYY(V)$. Its primitive part is $\Prim \YYY(V) = (\Prim_{As}\YYY)(V)$. It is a multibrace algebra for which  $M_{pq}=0$ for $p\geq 2$, therefore it is a brace algebra. Since these structures are functorial in $V$ we get a surjective map of operads
$$Brace \epi \Prim_{As}\YYY.$$
In order to show that this is an isomorphism, it suffices to show that, for any $n$, the spaces $Brace(n)$ and $(\Prim_{As}\YYY)(n)$ have the same dimension. Since we have isomorphisms $Dend \cong As^c \circ Brace$ and 
$\YYY \cong As^c \circ \Prim_{As}\YYY$, it suffices to show that 
$\dim Dend(n)= \dim \YYY(n) $. We claim that the isomorphism of CHAs $Dipt(V)\cong T^c(MB(V))\cong 2as(V)$ induces an isomorphism $Dend(V)\cong \YYY(V)$ which implies the expected equality. This last isomorphism is a consequence of the following Proposition.
\end{proo}

\begin{prop} The quotient of the operad $Dipt$ by the relations  $M_{pq}=0$ for $p\geq 2$ is the operad $Dend$.
\end{prop}
\begin{proo} We know that the operad $Dend$ is the quotient of $Dipt$ by the relation $M_{21}=0$ (cf.\ last line of section \ref{dipttoMB}). Therefore it suffices to show that, in $Dipt$, $M_{21}=0$ implies $M_{pq}=0$ for any $p\geq 2$. This is an immediate consequence of the following Lemma.
\end{proo}

\begin{lemma}\label{relationsdansDipt} In the free dipterous algebra $Dipt(V)$ the following relations hold:
$$
M_{n1}(x_{1}\ldots x_{n} ; y) = M_{(n-1)1}( (x_{1}\d x_{2})x_{3}\ldots x_{n} ; y) - x_{1}\d M_{(n-1)1}(x_{2}x_{3}\ldots x_{n} ; y)$$
for $n\geq 3$,
$$\displaylines{M_{nm}(x_{1}\ldots x_{n} ; y_{1}\ldots y_{m})= M_{n(m-1)} (x_{1}\ldots x_{n} ; (y_{1}\d y_{2})y_{3}\ldots y_{m})\hfill \cr
\hfill - y_{1}\d M_{n(m-1)} (x_{1}\ldots x_{n} ; y_{2})y_{3}\ldots y_{m})\cr
}$$
for $m>1$.
\end{lemma}
\begin{proo} We prove the first equality applying a recursive argument. The result  may be checked for $n=3$ by a straightforward calculation.

\noindent Let $n>3$ and suppose that the equality holds for all $M_{h1}$, with $h<n$. We get:
$$\displaylines {
\llap {(1)}\ \omega ^{\succ }(x_1\dots x_n)*y=
\sum _{r=0}^n\bigl (\sum _{j=0}^{n-r}\omega^{\succ }(x_1\dots x_rM_{j1}(\dots
x_{r+j};y)\dots x_n)\bigr ). \cr 
\llap {(2)}\ \omega ^{\succ }(x_1\dots x_n)* y=\hfill \cr
\omega ^{\succ }((x_1\succ x_2)\dots x_n)*y=\hfill \cr
\sum _{r=2}^n\bigl ( \sum _{j=0}^{n-r}\omega ^{\succ }(x_1\dots x_rM_{j1}(\dots x_{r+j};y)\dots x_n)\bigr )+\hfill\cr
\hfill \sum _{j=2}^n\omega ^{\succ }(M_{(j-1)1}((x_1\succ x_2)\dots x_j;y)\dots x_n)+
 \omega ^{\succ }(y(x_1\succ x_2)\dots x_n).\cr }$$

\noindent Equalities $(1)$ and $(2)$ imply:
$$\displaylines {
\llap (3)\  \sum _{j=0}^n\omega ^{\succ }(M_{j1}(x_1\dots x_j;y)\dots x_n)+
\sum _{j=0}^{n-1}\omega ^{\succ }(x_1M_{j1}(\dots x_{j+1};y)\dots x_n)=\hfill \cr 
\sum _{j=2}^n\omega ^{\succ }(M_{(j-1)1}((x_1\succ x_2)\dots x_j;y)\dots x_n)+ 
\omega ^{\succ }((y*x_1)x_2\dots x_n).\cr }$$ 
\medskip

\noindent The recursive hypothesis states that:
$$M_{j1}(x_1\dots x_j;y)+x_1\succ M_{(j-1)1}(x_2\dots x_j;y)=
M_{(j-1)1}((x_1\succ x_2)x_3\dots x_j;y),$$
for $3\leq j\leq n-1$.
\medskip

\noindent Applying it, $(3)$ becomes:
$$\displaylines {
M_{n1}(x_1\dots x_n;y)+ x_1\succ M_{(n-1)1}(x_2\dots x_n;y)+\omega ^{\succ}(yx_1\dots x_n)+\hfill \cr
\omega ^{\succ}(M_{11}(y;x_1)x_2\dots x_n)+\omega ^{\succ}(M_{11}(x_1\succ x_2;y)x_3\dots x_n)+
\omega ^{\succ}(x_1yx_2\dots x_n)=\hfill\cr
\hfill \omega ^{\succ}(M_{11}(x_1\succ x_2;y)x_3\dots x_n)+
M_{(n-1)1}((x_1\succ x_2)x_3\dots x_n;y)+
 \omega ^{\succ}((y*x_1)x_2\dots x_n).\cr }$$

\noindent But, since $M_{11}(y;x_1)=y* x_1-y\succ x_1 -x_1\succ y$, we get:
$$M_{n1}(x_1\dots x_n;y)+x_1\succ M_{(n-1)1}(x_2\dots x_n;y)=M_{(n-1)1}((x_1\succ x_2)x_3\dots x_n;y),$$
which ends the proof of the first equality.
\bigskip

For the second statement, the result is easy to verify for all $n\geq 0$ and $m=2$. 
For $m>2$, suppose that the equality holds for all $1\leq k\leq n$ and $1\leq h<m$.

We get:
$$\displaylines {
\llap {(4)}\ \omega ^{\succ}(x_1\dots x_n)*\omega ^{\succ}(y_1\dots y_m)=
\omega ^{\succ} (x_1\dots x_n)*\omega ^{\succ}((y_1\succ y_2)y_3\dots y_m)=\hfill \cr
\sum
\omega ^{\succ} (x_1\dots x_h M_{j(k-1)}(\dots x_{h+j};(y_1\succ y_2)
 \dots y_k) M_{r_1s_1}\dots M_{r_ls_l}(\dots  x_n \dots y_m))\cr
 \hfill +M_{n(m-1)}(x_{1}\dots x_{n};(y_{1}\succ y_{2})\dots y_{m}),\cr }$$
where the sum is taken over $0\leq h\leq n , 2\leq k\leq m , 0\leq j\leq n-h , j<n$ and the compositions $(r_{1}, \ldots, r_{l})$ of $n-h-j$, and $(s_{1}, \ldots, s_{l})$ of $m-k$. 
 Note that if $j=0$, then $k=2$.
 \medskip
 
\noindent Applying the formula recursively, we get:
\begin{enumerate}
\item if $h>0$, then 
$$\displaylines {
\omega ^{\succ }(x_{1}\dots x_{h}M_{j(k-1)}(\dots x_{h+j};(y_{1}\succ y_{2})\dots y_{k}))=\hfill\cr
\omega ^{\succ }(x_{1}\dots x_{h}M_{jk}(\dots x_{h+j};y_{1}\dots y_{k}))+
\omega ^{\succ}((\omega ^{\succ}(x_{1}\dots x_{h})*y_{1})M_{j(k-1)}(\dots x_{h+j};y_{2}\dots y_{k}))=\cr
\hfill \sum_{0\leq l\leq s\leq h} \omega^{\succ}(x_{1}\dots x_{l}M_{(s-l)1}(\dots x_{s};y_{1})\dots x_{h}M_{j(k-1)}(\dots x_{h+j};y_{2}\dots y_{k})).\cr}$$
\item if $h=0$ , then
$$\displaylines {
M_{j(k-1)}(x_{1}\dots x_{j};(y_{1}\succ y_{2})\dots y_{k})=\hfill\cr
\hfill M_{jk}(x_{1}\dots x_{j};y_{1}\dots y_{k})+y_{1}\succ M_{j(k-1)}(x_{1}\dots x_{j};y_{2}\dots y_{k}).\cr}$$
\end{enumerate}
\medskip

\noindent The formulas above imply the equality:
$$\displaylines {
\sum
\omega ^{\succ} (x_1\dots x_h M_{j(k-1)}(\dots x_{h+j};(y_1\succ y_2)
 \dots y_k) M_{r_1s_1}\dots M_{r_ls_l}(\dots  x_n \dots y_m))=\hfill\cr
\hfill \sum\omega^{\succ}(M_{i_{1}j_{1}}\dots M_{i_{r}j_{r}})(x_{1}\dots 
 x_{n}y_{1}\dots y_{m})- y_{1}\succ M_{n(m-1)}(x_{1}\dots x_{n};y_{2}\dots y_{m}).\cr }$$
where the first sum is taken over $0\leq h< n ,  2\leq k\leq m , 0\leq j\leq n-h$ and the compositions $(r_{1},\ldots , r_{l}), (s_{1},\ldots , s_{l})$, and the second sum is taken over $ k\geq 2 , i_{1}+\cdots +i_{r}=n , j_{1}+\cdots +j_{r}=m$.
\noindent Replacing in $(4)$, we get the equality:
$$\displaylines {
\sum _{k\geq 2}\bigl (\sum _{{i_{1}+\dots +i_{r}=n}\atop {j_{1}+\dots +j_{r}=m}}\omega^{\succ}(M_{i_{1}j_{1}}\dots M_{i_{r}j_{r}})(x_{1}\dots 
 x_{n}y_{1}\dots y_{m})+
 M_{nm}(x_{1}\dots x_{n};y_{1}\dots y_{m})=\hfill\cr
 \omega ^{\succ}(x_1\dots x_n)*\omega ^{\succ}(y_1\dots y_m)=\hfill\cr
\sum _{k\geq 2}\bigl (\sum _{{i_{1}+\dots +i_{r}=n}\atop {j_{1}+\dots +j_{r}=m}}\omega^{\succ}(M_{i_{1}j_{1}})\dots M_{i_{r}j_{r}})(x_{1}\dots 
 x_{n}y_{1}\dots y_{m})-\hfill\cr
 \hfill y_{1}\succ M_{n(m-1)}(x_{1}\dots x_{n};y_{2}\dots y_{m})+M_{n(m-1)}(x_{1}\dots x_{n};(y_{1}\succ y_{2})\dots y_{m}),\cr}$$
which implies:
$$\displaylines { M_{nm}(x_{1}\dots x_{n};y_{1}\dots y_{m})=\hfill\cr
\hfill y_{1}\succ M_{n(m-1)}(x_{1}\dots x_{n};y_{2}\dots y_{m})+M_{n(m-1)}(x_{1}\dots x_{n};(y_{1}\succ y_{2})\dots y_{m}),\cr}$$
as expected.
\end{proo}

\subsection{On the operad $\YYY$} As a corollary of Theorem \ref{tripleYYY} we deduce an isomorphism of $S$-modules:
$$\YYY \cong T^c(Brace).$$
Hence $\YYY(n)$ has the same dimension as $Dend(n)$, that is
$$\dim \YYY(n) = c_{n}\times n!.$$


 \section{(Co)free-(co)commutative CHA}\label{s:cofreecocom} 
 
    In this section we study the cofree-cocommutative  CHAs. We will see that they are $ComAs$-bialgebras. The space of primitives inherits a  symmetric multibrace (SMB) algebra structure. We omit the proofs which are completely analogous to the non-symmetric case. In this section $\KK$ is a characteristic zero field.

\subsection{Definition}\label{cofree-cocommutativeCHA} A \emph{cofree-cocommutative   combinatorial Hopf algebra} $\HH$ (cofree-cocommutative  CHA for short) is a cofree-cocommutative bialgebra with a prescribed isomorphism with $S^c(R)$ (here $R:=\Prim \HH$). Equivalently, it is a vector space $R$ with a product $*$ on $S^c(R)$ which makes $(S^c(R), *, \DD)$ ($\DD=$ deconcatenation) into a Hopf algebra.

\subsection{Symmetric multibraces} Let  $\HH=S^c(R)$ be a cofree-coassociative CHA. Since it is cofree, the multiplication $*: S^c(R)\t S^c(R)\to S^c(R)$ is completely determined by its projection onto $R$. So, there are well-defined maps $M_{pq}: S^p R\t S^q R\to R$, called \emph{symmetric multibrace operations}, given by the following composite:

$$\xymatrix{
S^c(R)\t S^c(R)\ar^-{*}[r] & S^c(R)\ar^{proj}[d] \\
  S^p R\t S^q R\ar^{inc}[u]\ar^-{M_{pq}}[r] & R
}$$  
  
In the next proposition we make explicit all the relations satisfied by the symmetric multibrace operations. We use the following notation. First $Sh(p,q)$ denotes the set of permutations made of $(p,q)$-shuffles. Second, 
for any integer $m$, a \emph{composition} of $m$ is an ordered sequence ${\underline k}=(k_1,\dots ,k_r)$ of nonnegative integers $k_i$ such that $\displaystyle {\sum _{i=1}^rk_i=m}$. We write $(k_1,\dots ,k_r)\vdash m$ when ${\underline k}$ is a composition of $m$. 
 
\begin{prop}\label{SRijk}  Let $M_{pq}: S^p R\t S^q R\to R$ be $(p,q)$-symmetric operations on the vector space $R$. Let $*$ be the unique binary operation on $S^cR$ which is a cocommutative coalgebra morphism and whose projection onto $R$ coincides with the $M_{pq}$'s. The following assertions are equivalent:

  (a) the operation $*$ is associative and unital,
  
  (b) the $(p,q)$-symmetric operations $M_{pq}$ satisfy the following relations $\SSS\RRR$:
  $$M_{00}=0, M_{01}= \id = M_{10}, M_{0q}= 0= M_{p0}, \textrm { for } p>1, q>1,$$
  and, for any integers $i,j,k$ greater than or equal to $1$ and  any elements\\
   $x_1,\dots ,x_n,y_1,\dots ,y_m,z_1,\dots ,z_r\in R$:
 
  $ (\SSS\RRR_{ijk}):$

$$\displaylines {
\sum _{r\geq 1}\bigl( \sum _{{(n_1,\dots ,n_k)\vdash n}\atop{(m_1,\dots ,m_k)\vdash m}}\frac{1}{k!}\bigl ( \sum _{{\sigma \in Sh(n_1,\dots ,n_k)}\atop {\tau\in Sh(m_1,\dots ,m_k)}}M_{kr}(M_{n_1m_1}
({\underline x}_1^{\sigma};{\underline y}_1^{\tau})\dots M_{n_km_k}({\underline x}_{k}^{\sigma};{\underline y}_{k}^{\tau});z_1,\dots ,z_r)\bigr )\bigr)=\hfill \cr
\sum _{l\geq 1}\bigl( \sum _{{(m_1,\dots ,m_l)\vdash m}\atop {(r_1,\dots ,r_l)\vdash l}}\frac{1}{l!}\bigl ( \sum _{{\gamma \in Sh(m_1,\dots ,m_l)}\atop {\delta\in Sh(r_1,\dots ,r_l)}}M_{nl}(x_1\dots x_n;M_{m_1r_1}({\underline y}_1^{\gamma};{\underline z}_{1}^{\delta}),\dots ,M_{m_lr_l}({\underline y}_{l}^{\gamma};{\underline z}_{l}^{\delta})\bigr )\bigr ),\cr }$$
where the first sum is taken over all pairs of compositions
 $(n_1,\dots ,n_k)\vdash n$ and $(m_1,\dots ,m_k)\vdash m$ such that if $n_i=0$ then $m_i\neq 0$, for $1\leq i\leq k$, while the second one is taken over all pairs of compositions $(m_1,\dots ,m_l)\vdash m$ and $(r_1,\dots ,r_l)\vdash r$ such that if $m_j=0$ then $r_j\neq 0$, for $1\leq j\leq l$. The elements 
${\underline x}_i^{\sigma}$, ${\underline y}_i^{\tau}$, ${\underline y}_j^{\gamma}$ and ${\underline z}_j^{\delta}$ are defined by:
$$w_i^{\omega}:=w_{\omega (p_1+\dots +p_{i-1}+1)}\t \cdots \t w_{\omega (p_1+\dots +p_i)},$$
for $w \in R^{\otimes p}$, $(p_1,\dots ,p_s)\vdash p$, $\omega \in {\mbox {\it Sh}(p_1,\dots ,p_s)}$, and $1\leq i\leq s$.
\end{prop}

  \begin{proo} From the unitality and counitality property of the Hopf algebra $S^c(R)$ we deduce that
$$M_{00}= 0, \ M_{10}= \id_V = M_{01}, \ {\rm and }\ M_{n0}=0=M_{0n}\ 
{\rm for} \ n\geq 2.$$
  
The operation $*$ is related to the operations $M_{pq}$ by the following formula:
$$u_{1}\ldots u_{p}* v_{1}\ldots v_{q}= \sum_{k=1}^{p+q}(u_{1}\ldots u_{p}* v_{1}\ldots v_{q})_{k}$$
where the component $(u_{1}\ldots u_{p}* v_{1}\ldots v_{q})_{k}\in S^k(R)$ is given by:
\begin{displaymath} \begin{array}{l}
(u_{1}\ldots u_{p}* v_{1}\ldots v_{q})_{k}=\\
\sum M_{i_1j_1}(u_1\ldots u_{i_1},v_1\ldots v_{j_1})M_{i_2j_2}(u_{i_1+1}\ldots u_{i_1+i_2},v_{j_1+1}\ldots 
v_{j_1+j_2})\ldots \qquad\qquad \\
\hfill \ldots M_{i_kj_k}( \ldots u_p, \ldots v_q)\in  \Vt k\ .
\end{array} \end{displaymath}
where the sum is extended to all the sequences of indices satisfying the following conditions:

-- each sequence of integers $\{1, \ldots i_{1}\} , \{i_{1}+1, \ldots, i_{1}+i_{2}\}, \ldots $ is ordered

-- $1\leq i_{1}+1\leq i_{1}+i_{2}+1\leq \ldots $,

-- each sequence of integers $\{1, \ldots j_{1}\} , \{j_{1}+1, \ldots, j_{1}+j_{2}\}, \ldots $ is ordered

-- $1\leq j_{1}+1\leq j_{1}+j_{2}+1\leq \ldots $.

This formula, which holds in the cocommutative case, is analogous to  the formula given in the coassociative case in section 1.4 of \cite{LRstr} (cf.\ Proposition \ref{Rijk}). The difference (restriction on the sequences of indices) is due to the fact that the coproduct is cocommutative, hence the map $\Delta:S(V)\to  S(V)\t S(V)$ is viewed as a composite $\Delta:S(V)\to  S^2(S(V)\t S(V))\mono S(V)\t S(V)$. 

For instance we get 
$$(u_{1}\ldots u_{p}* v_{1}\ldots v_{q})_{1}=M_{pq}(u_{1}\ldots u_{p} ; v_{1}\ldots v_{q}),$$
$$(u_{1}\ldots u_{p}* v_{1}\ldots v_{q})_{p+q}=u_{1}\ldots u_{p}v_{1}\ldots v_{q},$$
$$(u_{1}u_{2}* v_{1})_{2}=u_{1}M_{11}(u_{2};v_{1})+ M_{11}(u_{1}; v_{1})u_{2}.$$

By computing the component in $R$ of the two elements:
$$(x_1\cdots x_i * y_1 \cdots y_j)*z_1\cdots z_k = x_1\cdots x_i *( y_1 
\cdots y_j*z_1\cdots z_k)$$
we get the relation $ (\SSS\RRR_{ijk})$. The rest of the proof is as in Proposition \ref{Rijk}.
  \end{proo}
  
    \subsection{Example}\label{ExSR} The first nontrivial relation is $(\SSS\RRR_{111})$ which reads
  $$M_{21}(uv;w) + M_{11}(M_{11}(u;v); w) =  M_{11}(u; M_{11}(v;w))+ M_{12}(u; vw)$$
  for any $u,v,w\in R$. Remember that we are working in the symmetric framework, so $uv=vu$.
Here are the next ones:

\medskip

\noindent $(\SSS\RRR_{112})$
$$\displaylines{M_{22}(uv;wx)+M_{12}(M_{11}(u;v);wx)=\hfill \cr
\hfill M_{13}(u;vwx)+M_{12}(u;M_{11}(v;w)x)+M_{12}(u;wM_{11}(v;x))+M_{11}(u;M_{12}(v;wx)),}$$

\noindent $(\SSS\RRR_{121})$
$$\displaylines{M_{31}(uvw;x)+M_{21}(M_{11}(u;v)w;x)+M_{21}(vM_{11}(u;w);x)+M_{11}(M_{12}(u;vw);x)=\hfill \cr
\hfill M_{13}(u;vwx)+M_{12}(u;vM_{11}(w;x))
+M_{12}(u;M_{11}(v;x)w)+M_{11}(u;M_{21}(vw;x)),}$$

\noindent $(\SSS\RRR_{211})$
$$\displaylines{M_{31}(uvw;x)+M_{21}(uM_{11}(v;w);x)+M_{21}(M_{11}(u;wv);x)+M_{11}(M_{21}(uv;w);x)=\hfill \cr 
\hfill M_{22}(uv;wx)+M_{21}(uv;M_{11}(w;x)),}$$

    \subsection{Symmetric multibrace algebra}\label{SMB} By definition a \emph{symmetric multibrace algebra} (or $SMB$-algebra) is a vector space $R$ equipped with symmetric multibrace operations $M_{pq}$, that is $p+q$-ary operations $M_{pq}:S^p R\t S^q R\to R$, which satisfy all the conditions $\SSS\RRR$ of Proposition \ref{SRijk}.
 
   \begin{thm}\label{coComCombSMB} There is an equivalence of categories between the cofree-cocom\-mutative  CHAs and the SMB-algebras.
  \end{thm}
   
\subsection{Commutative-associative bialgebra}\label{ComAs} A \emph{$ComAs$-algebra} is a vector space equipped with two associative operations, one of them being commutative. We denote by $x\cdot y$ the commutative one and by $x*y$ the associative one. A $ComAs$-algebra is unital if there is an element $1$ which is a unit for both operations.

A \emph{$ComAs$-bialgebra} (i.e.\ a $Com^c\pt ComAs$-bialgebra) is a unital $ComAs$-algebra equipped with a counital cocommutative coproduct which satisfies the Hopf compatibility relation for both operations.

\begin{prop} Any cofree-cocommutative  CHA is a $ComAs$-bialgebra.
\end{prop}

\begin{proo} Let $\HH=S^c(V)$ be a cofree-cocommutative CHA. We equip $S^c(V)=S(V)$ with the standard commutative product, denoted $\cdot$, so that $(S(V), \cdot, \DD)$ is the classical polynomial algebra. Since $\HH$ is endowed with another associative product, denoted $*$, we have a $ComAs$-bialgebra $\HH = (S(V), \cdot, *, \DD)$.
\end{proo}

\begin{prop} The free $ComAs$ algebra over $V$ is a natural  $ComAs$-bialgebra.
\end{prop}

\begin{proo} The tensor product of two $ComAs$-algebras is still a $ComAs$-algebra. So we can apply the method of \cite{JLLsci} to construct the expected coproduct.
\end{proo}

\subsection{From $ComAs$-algebras to $SMB$-algebras}\label{ComAsSMB}
Let us recall that a conilpotent commutative and cocommutative Hopf algebra  $(\HH, \cdot ,\Delta )$ is isomorphic to $S^c({\rm Prim}(\HH))$. We choose the Eulerian idempotent $e^{(1)}:\HH\longrightarrow \HH$ (cf.\ \cite{Reutenauer, JLLeuler}), to construct this isomorphism.

If $\HH$ is a $Com^c\pt ComAs$-bialgebra, then, a fortiori, it is a $Com^c\pt Com$-bialgebra and we identify $\HH$ with $S^c({\rm Prim}(\HH))$ as above. We identify the dot product $x\cdot y$ with the polynomial product $xy$. The associative product $*$ induces a $SMB$ structure on ${\rm Prim}(\HH)$. As in the case of $2$-associative bialgebras, the relationship between the associative product, the commutative product and the multibraces is given by:
$$\displaylines { (x_1\cdots x_n) *(y_1\cdots y_m)=\hfill \cr
\hfill \sum _r\frac {1}{r!}\bigl(\sum _{{(k_1,\dots ,k_r)\vdash n}\atop {(h_1,\dots ,h_r)\vdash m}}\bigl (\sum _{{\sigma \in Sh(k_1,\dots ,k_r)}\atop
{\tau \in Sh(h_1,\dots ,h_r)}} M_{n_1m_1}({\underline x}_1^{\sigma},{\underline y}_1^{\tau})\cdots M_{n_rm_r}({\underline x}_r^{\sigma},{\underline y}_r^{\tau})\bigr)\bigr),\cr }$$
where the sum is taken over all pairs of compositions $(k_1,\dots ,k_r)\vdash n$ and $(h_1,\dots ,h_r)\vdash m$ such that if $k_i=0$ then $h_i\neq 0$.

Note that the formula for $M_{nm}$ may be obtained easily in terms of $*$, $\cdot $ and the $M_{ij}$\rq s, for $1\leq i\leq n$, $1\leq j\leq m$ and  $i+j<n+m$. For instance:
$$\displaylines{
M_{11}(x_1; y_1)=x_1*y_1-x_1y_1,\hfill \cr
M_{12}(x_1; y_1y_2)= x_1*(y_1y_2)-x_1y_1y_2-M_{11}(x_1; y_1)y_2-y_1M_{11}(x_1; y_2)\hfill \cr
\qquad = x*(yz) - (x*y)z - y(x*z) +xyz.\cr
}$$

\begin{prop} There is a forgetful functor $F$ from the category of $ComAs$-algebras to the category of $SMB$-algebras.
\end{prop}

\begin{thm}\label{tripleComAs} There is a good triple of operads
$$(Com, ComAs, SMB).$$
For any $V$ there is a canonical isomorphism
$$ComAs(V) \cong S^c(SMB(V)).$$
\end{thm}   

\subsection{The operad $ComAs$}\label{operadComAs} The operad $ComAs$ is a set-theoretic operad. Its underlying set was described by R.\ Stanley  in \cite{Stanley74}. The dimension of $ComAs(n)$ is the number $d_{n}$ of  ``labeled series-parallel posets with $n$ points'': 
$$\{1, 3, 19, 195, 2791, 51303, \ldots\quad \}.$$

\subsection{The operad $SMB$}\label{operadSMB} The operad $SMB$ can be described by using the ``connected labeled series-parallel posets with $n$ nodes'', cf.\ \cite{Stanley74}. Let us mention the first dimensions $f_{n}=\dim SMB(n)$ for $n\geq 1$:
$$\{1,2,12,122,1740, 31922, \ldots \quad \} .$$

From the presentation of the operad $SMB$ it is clear that we can modify it by reducing the number of generators. For instance the relation $\SSS\RRR_{111}$ shows that we can get rid of one of the two generators $M_{12}, M_{21}$. More symmetrically we can replace the two generators $M_{12}$ and $M_{21}$ by $M_{12}+M_{21}$. This discussion is postponed to a further paper.
  
  \subsection{From Lie algebras to $SMB$-algebras}\label{LieSMB}  Let ${\gg}$ be a Lie algebra and let $U({\gg})$ be its universal enveloping algebra. By the Cartier-Milnor-Moore theorem we know that $U({\gg})$ is isomorphic to $S^c(\gg)$ as a coalgebra. Let us choose the Eulerian idempotents to realize this isomorphism explicitly, cf.\ \cite{Reutenauer, JLLeuler}:
  $$e:U(\gg) \cong S^c(\gg).$$
  Recall from \cite{JLLeuler} that the component $e^{(1)}$, whose image lies in $\gg= \Prim U(\gg) \subset U(\gg)$ is obtained by 
  $$e^{(1)}= \log^{\star}(\Id),$$
  where $ \log^{\star}$ is the convolution log. The other components are given by 
  $$e^{(i)}= \frac{(e^{(1)})^{\star i}}{i!}.$$
  Once this isomorphism is chosen, Theorem \ref{coComCombSMB} implies that  there is an SMB-algebra structure on $\gg$. Hence we have constructed a forgetful functor
  $$Lie\alg \to SMB\alg.$$
  In low dimension we get the following formulas:
  $$M_{11}(x ; y) = \frac{1}{2}[x,y],$$
  $$M_{12}(x ; yz) = \frac{1}{6}[[x,y],z] -  \frac{1}{12}[x,[y,z]],$$
    $$M_{21}(xy ; z) = -\frac{1}{12}[[x,y],z] +\frac{1}{6}[x,[y,z]].$$
  
\subsection{Free-commutative CHA}\label{freecomCHA} As in the preceding examples there is a dual version of the results of this section. We leave it to the interested reader to phrase them in details.

\section{(Co)free-(co)commutative right-sided CHA}\label{s:cofreecocomrs} 

 In this section we study the cofree-cocommutative  CHAs which satisfy the right-sided property {\bf (r-s)}.  The space of primitives inherits a symmetric brace algebra structure, which turns out to be the same as a pre-Lie algebra structure. As shown in \cite{ChapotonLivernet01} the free pre-Lie algebra on one generator gives rise to a cofree-cocommutative CHA which is the Grossman-Larson algebra \cite{GrossmanLarson}. In the dual framework, the cofree pre-Lie coalgebra on one generator gives rise to a free-commutative CHA which is the Connes-Kreimer Hopf algebra \cite{CK}.  In this section $\KK$ is a characteristic zero field.  Examples will be given in section  Examples will be given in section \ref{Exs}. They include the dual of the Fa\`a di Bruno algebra (\ref{FaadiBruno}), the symmetric functions algebra (\ref{sym}).

\subsection{Right-sided cofree-cocommutative  CHA}\label{rsccCHA} Let $\HH$ be a cofree-cocom\-mutative  CHA. By definition $\HH=(S^c(R), *, \DD)$ is said to be \emph{right-sided}  if the following condition holds:
 
 \medskip
 
 \noindent {\bf (r-s)} for any integer $q$ the subspace $\bigoplus_{n\geq q} {\HH^n}$ is a right-sided ideal of $\HH$,
 
 \medskip
 
\noindent  where  $\HH^n:= S^n(R)$. 
We have seen that the associative product $*$ is given by a family 
$\{M_{pq}\}_{p,q\geq 0}$ of  symmetric multibrace operations (cf.\ Proposition  \ref{SRijk}). 
If $\HH$ is right-sided, then, as in \ref{rsbrace}, the operations $M_{pq}$ are 0 for any $p\geq 2$. Hence we get a \emph{symmetric brace algebra}, $SBrace$-algebra for short, that is a vector space $R$ equipped with operations $M_{1q}$ verifying:

$$\displaylines {
M_{1m}(M_{1n}(x;y_1\dots y_n);z_1\dots z_m)=\hfill\cr
\sum _{(m_1,\dots ,m_{n+1}\vdash m}\bigl ( \sum _{\sigma \in Sh(m_1,\dots ,m_{n+1})}M_{1(n+m_{n+1})}(x;M_{1m_1}(y_1,{\underline z}_1^{\sigma })\dots 
M_{1m_n}(y_n;{\underline z}_n^{\sigma }){\underline z}_{n+1}^{\sigma})\bigr ),\hfill\cr }$$
where ${\underline z}_i^{\sigma }:=z_{\sigma (m_1+\dots +m_{i-1}+1)}\otimes \dots
\otimes z_{\sigma (m_1+\dots +m_i)}$.

We remark that the first relation ($n=1$) implies that the operation $M_{11}$ is pre-Lie. Recall that an operation $\{-,-\}$ is \emph{pre-Lie} if
  $$\{\{x,y\},z\}-\{x,\{y,z\}\}= \{\{x,z\},y\}-\{x,\{z,y\}\}  .$$
  Moreover it follows from $\SSS\RRR_{111}$ (cf.\ \ref{ExSR}) that the operation $M_{12}$ is completely determined by the operation $M_{11}$ since $M_{21}=0$. More generally all the operations $M_{1n}$ are determined by $M_{11}$ and the following result holds.

\begin{prop} [Guin-Oudom \cite{GuinOudom05, GuinOudom08}]\label{Sbrace=preLie} A symmetric brace algebra is equivalent to a pre-Lie algebra.
\end{prop}

\begin{proo} A proof can be found in \cite{GuinOudom08}, see also \cite{LadaMarkl}. The operations $M_{1n}$ are obtained from $M_{11}$ recursively by the formulas:
$$\begin{matrix}
M_{11}(x;y)&:=&\{ x,y\},\\
M_{1n}(x;y_1\dots y_n)&:=&M_{11}(M_{1(n-1)}(x;y_1\dots y_{n-1});y_n)-\\
&&{\displaystyle \sum _{1\leq i\leq n-1}}M_{1(n-1)}(x;y_1\dots M_{11}(y_i;y_n)\dots y_{n-1}),\end{matrix}$$
for $x,y_1,\dots ,y_n\in V$.
\end{proo}

 \begin{thm}\label{thm:cocomCHApreLie} The category of cofree-cocommutative right-sided CHAs is equivalent to the category of pre-Lie algebras.
\end{thm}
\begin{proo} Taking the primitives gives a functor from cofree-cocommutative r-s CHAs to pre-Lie algebras. Since any pre-Lie algebra gives rise to a symmetric brace algebra, we can apply the reconstruction functor to get a cofree-cocommutative CHA (cf.\ \ref{coComCombSMB}). Since we started with a symmetric brace algebra rather than a general SMB-algebra, the Hopf algebra is right-sided, hence the existence of a functor from pre-Lie algebras to cofree-cocommutative r-s CHAs. Obviously these two functors are inverse to each other.
\end{proo}

\subsection{Unreduced trees and free pre-Lie algebra}\label{unreducedtrees} In \ref{uptrees} we defined the planar unreduced trees. The quotient by the obvious relation gives the notion of (non-planar) \emph{unreduced trees}. For instance, in this setting, the two drawings

 \unreducedtreesdeux

\noindent define the same element. The set of unreduced trees with $n$ vertices is denoted $UT_{n}$. The free pre-Lie algebra on $V$ has been described by Livernet and Chapoton in \cite{ChapotonLivernet01} in terms of unreduced trees labeled by elements of $V$:
$$PreLie(V)= \sum_{n}\KK[UT_{n}]\t V^{\t n}.$$
The pre-Lie product is given by the following rule. If $\omega$ and $\omega '$ are two labelled unreduced trees, then the pre-Lie product $\{\omega,\omega '\}$ is the sum of all the trees obtained from  $\omega$ and $\omega '$ by drawing an edge from a vertex of $\omega$ to the root of $\omega '$. The root of this new tree is the root of $\omega$. For instance:

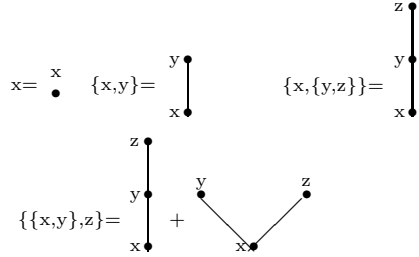
\begin{picture}(100,50)

\put(00,10) {\scriptsize{x=}}
\put(15,07){\scriptsize$\bullet$}
\put(15,15){\scriptsize x}

\put(30,10) {\scriptsize{\{x,y\}=}}
\put(65,00){\scriptsize$\bullet$}
\put(60,00){\scriptsize x}
\put(65,20){\scriptsize $\bullet$}
\put(60,20){\scriptsize y}
\put(67,00){\line(0,1){20}}

\put(100,10){\scriptsize{ \{x,\{y,z\}\}=}}
\put(150,00){\scriptsize $\bullet$}
\put(145,00){\scriptsize x}
\put(150,20){\scriptsize $\bullet$}
\put(145,20){\scriptsize y}
\put(150,40){\scriptsize $\bullet $}
\put(145,40){\scriptsize z}
\put(152,00){\line(0,1){20}}
\put(152,20){\line(0,1){20}}
\end{picture}

\begin{picture}(80,50)
\put(00,10){\scriptsize{ \{\{x,y\},z\}=}}

\put(50,00){\scriptsize $\bullet$}
\put(45,00){\scriptsize x}
\put(50,20){\scriptsize $\bullet$}
\put(45,20){\scriptsize y}
\put(50,40){\scriptsize $\bullet $}
\put(45,40){\scriptsize z}
\put(52,00){\line(0,1){20}}
\put(52,20){\line(0,1){20}}

\put(60,10){\scriptsize{+}}

\put(90,00){\scriptsize $\bullet$}
\put(85,00){\scriptsize x}
\put(70,20){\scriptsize $\bullet$}
\put(70,25){\scriptsize y}
\put(110,20){\scriptsize $\bullet$}
\put(110,25){\scriptsize z}
\put(92,00){\line(-1,1){20}}
\put(90,00){\line(1,1){20}}
\end{picture}

\subsection{The conjectural operad $\XXX$} The results of the three other cases lead us naturally to conjecture the existence of an operad structure $\XXX$ on $Com\circ preLie$ compatible with the operad structure  of $Com$ and $preLie$ (extension of operads, see 3.4.2 in \cite{JLLgbto}). Observe that

$$\dim \XXX(n) = (n+1)^{n-1}.$$

If so, then we conjecture that any cofree-cocommutative right-sided CHA is a $\XXX$-algebra.

 \subsection{Grossman-Larson Hopf algebra} In \cite{GrossmanLarson} Grossman and Larson constructed a cocommutative Hopf algebra on some trees. It turns out that this is exactly the combinatorial Hopf algebra $S^c(preLie(\KK))$.
 
 \subsection{Free-commutative right-sided CHAs and Connes-Kreimer algebra}\label{freecomrsCHA} As before the results of this section can be linearly dualized. It should be said that many examples in the literature appear as free-commutative (rather than cofree-cocommutative). For instance the dual of $S^c(preLie(\KK))$, that is $S(copreLie(\KK))$  is the Connes-Kreimer Hopf algebra (cf.\ \cite{Panaite, Hoffman00}). It can be constructed either directly as in \cite{CK}, or by means of the cofree pre-Lie coalgebra on one generator. 
 
 The direct construction consists in taking the free-commutative algebra over the set of unreduced trees. The coproduct is given by the formula:
 $$\DD(t)=\sum_{c}P_{c}(t) \t R_{c}(t),$$
where $P_{c}(t) $ is a polynomial of trees and $R_{c}(t) $ is a tree, the sum is over admissible cuts $c$. Let us recall that a cut of the tree $t$ is admissible if there is one and only one cut on any path from the root to a leaf. Among the pieces, one of them contains the root, this is the tree $R_{c}(t)$. The other ones assemble to give a polynomial, which is $P_{c}(t)$. Observe that there are two extreme cuts: under the root (it gives the element $t\t 1$), above the leaves (it gives the element $1\t t$).
 
 We end this section with the dual version of Theorem \ref{thm:cocomCHApreLie}.

 \begin{thm}\label{thm:comCHApreLie} The category of free-commutative right-sided CHAs is equivalent to the category of pre-Lie coalgebras.
\end{thm}
 
 \subsection{Cofree-cocommutative CHAs being right and left-sided} Let $\HH$ be a cofree-cocommutative CHA which is both right-sided and left-sided. As a consequence its primitive part is such that $M_{nm}=0$ except possibly for $M_{11}$. By formula $\RRR_{111}$ it follows that $M_{11}$ is associative.
 
 This case has been studied in \cite{JLLctd} where it is proved that there is a good triple of operads
 $$(Com, CTD, As)$$
 where $As$ is the operad of nonunital associative algebras, and $CTD$ is the operad of \emph{commutative tridendriform algebras}. Recall from \cite{LRtri} that a tridendriform algebra is determined by three binary operations $x\g y, x\d y, x\cdot y$ satisfying 7 relations. It is said to be \emph{commutative} whenever $x\g y = y\d x$ and $x\cdot y = y\cdot x$ for any $x$ and $y$. Then, the 7 relations come down to the following 4 relations:
 $$\displaylines{
 (x\g y)\g z = x\g (y\g z + z\g y + y\cdot z), \cr
  (x\g y)\cdot z = x\cdot  (z\g y), \cr
   (x\cdot y)\g z = x\cdot (y\g z), \cr
     (x\cdot y)\cdot z = x\cdot (y\cdot z), \cr 
 }$$
 The second relation has been, unfortunately, mistakenly omitted in \cite{JLLctd}.
 
 There is an equivalence of categories between the right-and-left-sided cofree-cocommutative CHAs and the nonunital associative algebras.
 
\section{Examples}\label{Exs}
  
There are many examples of combinatorial Hopf algebras in the literature. We already described some of them: the free algebras of type $Dipt$, $2as$, $ComAs$, $Dend$, $\YYY$ and their dual. In these cases the combinatorial objects at work are trees (of various types). The examples coming from quantum field theory are usually free-associative or free-commutative (for instance Connes-Kreimer algebra). In this section we list some specific examples and refer to the literature for more information. In \ref{MalvenutoReutenauer} we give an example of two different CHA structures on a given cofree-coassociative Hopf algebra.

 \subsection{Fa\`a di Bruno algebra}\label{FaadiBruno} On  $F=\bigoplus_{n\ge 1}\KK x_{n}$ the operation:
  $$\{x_{p}, x_{q}\} := -px_{p+q}$$
  is a pre-Lie product since
  $$\{\{x_{p},x_{q}\},x_{r}\} - \{x_{p},\{x_{q},x_{r}\}\}= p^2x_{p+q+r},$$
is symmetric in $q$ and $r$.

Up to a minor modification, this is the pre-Lie algebra coming from the nonsymmetric operad $As$, cf. \ref{operads}. 
 By Theorem \ref{thm:cocomCHApreLie} the coalgebra $S^c(F)$ is a combinatorial Hopf algebra which is cofree-commutative and right-sided as an algebra. This is exactly the graded dual of the Fa\`a di Bruno Hopf algebra, since the Lie bracket is given by $[x_{p},x_{q}]=(-p+q)x_{p+q}$, see for instance \cite{JoniRota, FGB, FGBV}. More precisely, let us introduce the dual basis $a_{n}:= \frac{1}{n!}x^*_{n-1}$. The linear dual $S^c(F)^*$ is the polynomial algebra in $a_{2}, \ldots, a_{n}, \ldots $ and the coproduct is given by 
 $$\DD(a_{n}) = \sum_{k=1}^n \sum_{\lambda} \binom{n}{\lambda ; k} a_{1}^{\lambda_{1}}\cdots a_{n}^{\lambda_{n}} \t a_{k}$$
 where $ \binom{n}{\lambda ; k}= \frac{n!}{\lambda_{1}! \cdots \lambda_{n}! (1!)^{\lambda_{1}}\cdots  (n!)^{\lambda_{n}}}$ for $\lambda_{1}+\cdots +\lambda_{n}= k$ and $\lambda_{1}+2\lambda_{2}+\cdots +n\lambda_{n}= n$.
 We see immediately on this formula that the Fa\`a di Bruno Hopf algebra is a right-sided coalgebra.
 
  This Hopf algebra proves helpful in studying the higher derivatives of a composition of formal power series (Fa\`a di Bruno's formula) because it is the Hopf algebra of functions on the group of invertible formal power series in one variable.
 
 It is immediate to check that the pre-Lie algebra $F$ is generated (in characteristic zero) by one element, namely $x_{1}$. Hence $F$ is a quotient of $PreLie(\KK)$. As a consequence the Fa\`a di Bruno Hopf algebra gets identified with a subalgebra of the Connes-Kreimer algebra.
 
 There is an analogue of the Fa\`a di Bruno algebra in the noncommutative framework due Brouder, Frabetti and Menous, cf.\ \cite{BFM}. Since $F$ is a pre-Lie algebra, a fortiori it is a brace algebra. Therefore, taking the cofree-coassociative coalgebra gives a combinatorial Hopf algebra. The graded dual is a noncommutative analogue of the Fa\`a di Bruno algebra.
 
\subsection{The algebra of quasi-symmetric functions ${\mbox {\bf QSym}}$}\label{qsym} (see \cite{HNT08}). 
Let $\KK [\xi_1,\dots ,\xi_n,\dots ]$ be the algebra of all polynomial functions on an infinite family of variables 
$\{\xi _n\}_{n\geq 1}$. The algebra of \emph{quasi-symmetric functions} ${\mbox {\bf QSym}}$ is the subalgebra of all polynomial functions $f\in \KK [\xi_1,\dots ,\xi_n,\dots ]$ such that 
$f=\sum _{n_1\dots n_r} f_{n_1\dots n_r}\bigl (\sum _{i_1<\dots <i_r}\xi_{i_1}^{n_1}\dots \xi_{i_r}^{n_r}\bigr )$, for certain coefficients $f_{i_1\dots i_r}\in \KK$. The product $*$ of ${\mbox {\bf QSym}}$ is the usual product of polynomials.

The subspace of homogeneous elements of degree $n$ of ${\mbox {\bf QSym}}$ has a natural basis $\{ x_{n_1\dots n_r}\}_{n_1+\dots +n_r=n}$, with:
$$x_{n_1\dots n_r}:=\sum _{i_1<\dots <i_r}\xi_{i_1}^{n_1}\dots \xi_{i_r}^{n_r}.$$

The coproduct on ${\mbox {\bf QSym}}$ is defined by:
$$\Delta (x_{n_1\dots n_r}):=\sum _{0\leq j\leq r}x_{n_1\dots n_j}\otimes x_{n_{j+1}\dots n_r}.$$

It is clear that the subspace of primitive elements of ${\mbox {\bf QSym}}$ is the vector space $\bigoplus _{n\geq 1}\KK x_n$. Moreover, the coalgebra isomorphism with
$T^c(\bigoplus _{n\geq 1}\KK x_n)$ sends the element $x_{n_1\dots n_r}$ to $x_{n_1}\otimes \dots 
\otimes x_{n_r}$.
\medskip

Note that ${\mbox {\bf QSym}}$ is cofree-coassociative and its product is associative and commutative. The multibrace structure of $\bigoplus _{n\geq 1}\KK x_n$ associated to the product $*$ is given by:
$$M_{nm}(x_{i_1}\dots x_{i_n};x_{j_1}\dots x_{j_m})=\begin{cases}x_{i_1+j_1},&{\rm for}\ n=m=1,\\
0,&{\rm otherwise.}\end{cases}$$
\medskip

The graded dual of ${\mbox {\bf QSym}}$ is the Solomon descent algebra ${\mbox {\bf NSym}}$.
It is the free algebra over the space $\bigoplus _{n\geq 1}\KK x_n$, with the coassociative cocommutative coproduct $\Delta ^*$ given by:
$$\Delta ^*(x_n):=\sum _{i=0}^nx_i\otimes x_{n-i}.$$

The Solomon descent algebra is a cofree-cocommutative  CHA, isomorphic as a coalgebra to $S^c\big(Lie(\bigoplus _{n\geq 1}\KK x_n)\big)$. This is a particular example of the case discussed in \ref{LieSMB}.

\subsection{The Hopf algebra ${\mbox {\bf Sym}}$}\label{sym} The Hopf algebra ${\mbox {\bf Sym}}$ is a cofree-cocommutative  CHA isomorphic as a coalgebra to $S^c\big(\bigoplus _{n\geq 1}\KK x_n\big)$. This is a particular example of the case discussed in \ref{LieSMB}.

\subsection{The Malvenuto-Reutenauer algebra $\HH_{MR}$}\label{MalvenutoReutenauer}(See \cite{MR}). The underlying space of $\HH_{MR}$ is the graded space $\bigoplus _{n\geq 0}\KK [S_n]$, where $S_n$ denotes the group of permutations of $n$ elements. The product of $\HH_{MR}$ is the shuffle product, given by:
$$\sigma *\tau :=\sum _{\delta \in Sh(n,m)}(\sigma \times \tau)\cdot \delta ^{-1},$$
for $\sigma \in S_n$ and $\tau \in S_m$; where ${\mbox {\it Sh}(n,m)}$ denotes the set of all ${\mbox {\it (n,m)}}$-shuffles and $\times $ denotes the concatenation of permutations.

Given a permutation $\sigma \in S_n$ and an integer $0\leq i\leq n$, let $\sigma _{(1)}^i$ be the map from $\{ 1,\dots ,i\}$ to $\{1,\dots ,n\}$ whose image is $\sigma _{(1)}^i=(\sigma (1),
\dots ,\sigma (i))$ and $\sigma _{(2)}^i$ be the map whose image is $\sigma _{(2)}^i=(\sigma (i+1),\dots ,\sigma (i))$.

The coproduct is defined as follows:
$$\Delta (\sigma):=\sum _{i=0}^n {\mbox{\it std}}(\sigma _{(1)}^i)\otimes {\mbox{\it std}}(\sigma _{(2)}^i),$$
where ${\mbox{\it std}}(\tau)$ is the unique permutation in $S_n$ such that ${\mbox{\it std}}(\tau)(i)<{\mbox{\it std}}(\tau)(j)$ if, and only if $\tau (i)<\tau (j)$, for any injective map $\tau :\{1,\dots ,k\}\longrightarrow \{1,\dots ,n\}$.
In degree two the primitive space is 1-dimensional generated by $12-21$. In degree three it is 3-dimensional generated by the elements:

\xymatrix@R=1pt @C=1pt{
u:=213-312,\\
v:=231 -132,\\
w:=321 - 132 - 213 + 123.}

An element $\sigma \in S_n$ is called {\it irreducible} (or connected) if $\sigma \notin \bigcup _{1\leq i\leq n-1}S_i\times S_{n-i}$. We denote by ${\mbox {\it Irr}}_n$ the subset of irreducible elements of $S_n$. In low dimension we get $Irr_{2}=\{21\}, Irr_{3}=\{231,312,321\}$.

We will give two different structures of combinatorial Hopf algebra on $\HH_{MR}$. For the first one it appears as a cofree-coassociative general CHA and for the second one it appears as a  cofree-coassociative right-sided CHA.

\noindent 1) Define an isomorphism $\varphi$ 
from $T^c\big(\bigoplus _{n}\KK[{\mbox {\it Irr}}_n]\big)$ to $\HH_{MR}$ by:
$$\varphi (x)=\begin{cases}e_{inf}^{(1)}(\sigma ) ,&{\rm for}\ x=\sigma \in \bigcup _n{\mbox {\it Irr}}_n\\
e_{inf}^{(1)}(\sigma _1)\times \dots \times e_{inf}^{(1)}(\sigma _m),&{\rm for}\ x=\sigma _1\otimes \cdots \t \sigma _m,\end{cases}$$
where 
$$e_{inf}^{(1)}(\sigma):=\sum _{i\geq 1}(-1)^{i+1}\times ^i\circ {\overline {\Delta }}^i(\sigma)= \sigma - \sigma_{(1)}\times \sigma_{(2)}+ \cdots $$
 is the (infinitesimal) idempotent defined in \cite{LRstr}. For instance, in low dimension, we get:

\xymatrix@R=1pt @C=1pt{
\varphi(21)= 21-12,\\
\varphi(231)=231-132-123+123= 231-132,\\
\varphi(312)=312-123-213+123=312-213,\\
\varphi(321)= 321-132-213+123.
}

Under this construction $(\HH_{MR}, \phi)$ become a cofree-coassociative general CHA. 
Hence the space $\bigoplus _{n\geq 1}\KK[{\mbox {\it Irr}}_n]$ inherits a structure of $MB$-algebra (see \cite{Ronco07} for explicit formulas). The  $2as$-bialgebra structure is given by the products $\times $ and $*$, and the coproduct $\Delta$.

\medskip

\noindent 2) Define an isomorphism $\theta : T^c\big(\bigoplus _{n}\KK[{\mbox {\it Irr}}_n]\big)\longrightarrow \HH_{MR}$ as follows:
$$\theta (x)=\begin{cases}e(\sigma ) ,&{\rm for}\ x=\sigma \in \bigcup _n{\mbox {\it Irr}}_n\\
((e(\sigma _1)\succ e(\sigma _2))\dots )\succ  e(\sigma _m),&{\rm for}\ x=\sigma _1\otimes \dots \t \sigma _m,\end{cases}$$
where $e(\sigma):=\sum _{i\geq 1}(-1)^{i+1}\omega _{\succ}^i\circ {\overline {\Delta }}^i(\sigma)$ is the idempotent defined in section \ref{pbtprimitives}.  So $(\HH_{MR}, \theta)$ is a cofree-coassociative right-sided CHA. Indeed, one can check that the multi-brace operations satisfy $M_{pq}=0$ for $p\geq 2$, so there is a brace structure on $\KK[{\mbox {\it Irr}}_n]$. The dendriform algebra structure on $\HH_{MR}$ is given by

$$\displaylines{
\sigma \succ \tau :=\sum _{{\delta \in Sh(n,m)}\atop {\delta (n+m)=n+m}}(\sigma \times \tau)\cdot \delta ^{-1},\cr
\sigma \prec \tau:=\sum _{{\delta \in Sh(n,m)}\atop {\delta (n)=n+m}}(\sigma \times \tau)\cdot \delta ^{-1}.\cr}$$

Observe that, because of the second condition under the summation sign, the sum is only over half-shuffles. In low dimension ($n\leq 3$) $\theta$ coincides with $\varphi$, but for $n\geq 4$ it is different.

The Hopf algebra of Solomon-Tits and the algebra ${\mbox {\bf FPQSym}}$ of parking functions may be described as $CHA$ in a similar way, cf.\ \cite{NovelliThibon, PalaciosRonco}.

In \cite{PoirierReutenauer} the authors show that there exists a sub-algebra of $\HH_{MR}$ which is spanned by the image of the standard Young tableaux. They construct a coalgebra structure for which the ``connected tableaux'' form a basis of the primitive part. Hence there is a $MB$-algebra structure on this latter space.

\subsection{Tensor module as indecomposables}\label{tensormodule}  In the literature there are several examples of combinatorial Hopf algebras whose space of indecomposables is the tensor module.

-- in \cite{BrouderSchmitt} Brouder and Schmitt construct a Hopf algebra structure on $T(\To(H))$, $S(\To(H))$ and on $S(\So(H))$ where $H$ is a non-unital bialgebra. These constructions generalize a construction of G.\ Pinter related to renormalisation in perturbative quantum field theory.

-- In \cite{Turaev} Turaev constructs a Hopf algebra structure on $T(T(V))$ (resp. $T(\To (V))$) which is free-associative right-sided. As we know there is a structure of pre-Lie coalgebra on $T(V)$ (resp. $\To(V)$) which is studied in details. It is part of a more general structure: a brace coalgebra structure, not studied in loc.cit.

\subsection{Quantum field theory} The Connes-Kreimer algebra is an example of the combinatorial Hopf algebras which appear in quantum field theory. They are based on graphs (Feynman graphs) and the rule is always the same: a subgraph is singled out and put on the right side. What is left is used to construct the left-hand side of the coproduct. See for instance \cite{ConnesMoscovici, FGB, BrouderFrabetti, Brouder, Tanasa}.

\subsection{Operads}\label{operads} For any nonsymmetric operad $\PP$, the space $\bigoplus_{n} \PP_{n}$ has a natural structure of brace algebra, therefore $T^c(\bigoplus_{n} \PP_{n})$ is a cofree-coassociative r-s CHA  (cf.\ for instance \cite{Holtkamp05, Moerdijk}). Similarly, if $\PP$ is a properad, then it can be shown that $\oplus_{n,m}\PP(n,m)$ is naturally equipped with a structure of $MB$-algebra, cf.\ \cite{MerkulovVallette}.  Therefore $T^c\big(\bigoplus_{n,m} \PP(n,m)\big)$ is a cofree-coassociative CHA.

\subsection{Incidence Hopf algebras}\label{incidence} Continuing the paper of Joni and Rota \cite{JoniRota} W.\ Schmitt studied in \cite{Schmitt94} the notion of incidence Hopf algebras. In the last sections of the work he described incidence Hopf algebras related to families of graphs closed under formation of induced subgraphs and sums. His examples include the the Fa\`a di Bruno algebra. The incidence Hopf algebras on graphs defined in \cite{Schmitt94} are either free-commutative or free-associative, so they enter in the 
examples studied in this paper.
\medskip

Consider for instance the free associative algebra spanned by the set ${\tilde {\mathcal G}}_0$ of all the isomorphism classes of connected simple graphs (i.e.\ having neither
 loops nor multiple edges) with linearly ordered vertex set, and let $H_l({\mathcal G})$ be the free associative algebra generated by ${\tilde {\mathcal G}}_0$. Note 
 that as a vector space $H_l({\mathcal G})$ is spanned by the isomorphism classes of simple graphs with linearly ordered vertex set.
 
\noindent For any simple connected graph $G$ with the set of vertices $V(G)$ linearly ordered and any subset of $U\subseteq V(G)$ the induced graph $G\vert U$ is the 
graph whose set of vertices is $U$, with the edge set formed by all the edges of $G$ which have both end-vertices in $U$. Note that $G\vert U$ is simple, 
and that $U$ inherits the linear order of $V(G)$. So, the isomorphism class of $G\vert U$ belongs to $H_l({\mathcal G})$. The coalgebra structure on $H_l({\mathcal G})$ is given by:
$$\Delta (\langle G\rangle):=\sum _{U\subseteq V(G)}\langle G\vert U\rangle\otimes \langle G\vert (V(G)- U)\rangle,$$
where $\langle-,-\rangle$ denotes the isomorphism class.

Clearly $H_l({\mathcal G})$ is free-associative and is cofree-cocommutative. Its graded dual $H_l({\mathcal G})^*$ is cofree-coassociative isomorphic to 
$T^C(\KK[{\tilde {\mathcal G}}_0])$. The multibrace structure on $\KK[{\tilde {\mathcal G}}_0]$ is given by:

\noindent $M_{nm}\big(\langle G_1\rangle\dots \langle G_n\rangle;\langle H_1\rangle\dots \langle H_m\rangle\big)=$ the sum of the isomorphism classes of all simple connected graphs $G$ such that:\begin{enumerate}
\item the set of vertices $V(G)$ is the disjoint union $\big(\bigcup _{i=1}^nV(G_i)\big)\bigcup \big(\bigcup _{j=1}^mV(H_j)\big)$,
\item If $i\neq k$, there does not exist any edge of $G$ between a vertex of $G_i$ and a vertex of $G_k$,
\item If $j\neq l$, there does not exist any edge of $G$ between a vertex of $H_j$ and a vertex of $H_l$,
\item the linear order on $V(G)$ is such that for $1\leq i<k\leq n$, the minimal vertex of $G_i$ is smaller than the minimal vertex of $G_k$; and for $1\leq j<k\leq m$, the minimal vertex of $H_j$ is smaller than the minimal vertex of $H_l$. Moreover, the induced orders on $G\vert G_i$ and $G\vert H_j$ coincide with the orders of $G_i$ and $H_j$, respectively.
\end{enumerate}

\subsection{Relationship with algebraic $K$-theory} The free dendriform algebra on one generator is spanned by the planar binary trees (cf.\ \ref{pbtDend} and \cite{JLLdig}). If, in each dimension, we take the sum of all the trees, then it spans a one-dimensional space, and the sum over all dimensions form a sub-Hopf algebra of the Loday-Ronco algebra $Dend(\KK)$. The same procedure applied to the free dendriform algebra over some decoration set $X$ provides a combinatorial Hopf algebra which is dual to the Hopf algebra constructed by Gangl, Goncharov and Levin in \cite{GGL}. The importance of this combinatorial Hopf algebra lies in its close relationship with computation in algebraic $K$-theory (cf.\ loc.cit.).

\section{Variations} Another interesting case consists in assuming both right-sidedness and left-sidedness on a CHA. It is treated in the first paragraph of this section.

 In this paper our object of study is conilpotent Hopf algebras. As we know it involves an associative product and a coassociative coproduct. But there exist more general types of bialgebras involving various kinds of operads. We briefly mention the pattern of a similar theory to generalized bialgebras in the second paragraph.

\subsection{Cofree-associative right and left-sided CHAs} Let us assume that the CHA $\HH$ is both right-sided and left-sided. Then it turns out that the multibrace operations are all trivial with the exception of $M_{11}$. From the formula \ref{FormulaforM11} we deduce that $M_{11}$ is an associative operation. So the primitive part is simply an associative algebra. This example is well-documented in the literature, it is called ``quasi-shuffle algebra'', cf.\ \cite{KEFGuo, JLLctd}.

\subsection{Generalized bialgebras} In \cite{JLLgbto} we introduced the notion of \emph{generalized bialgebras}, more specifically $\CCC^c\textrm{-}\AA$-bialgebras, where $\CCC$ governs the coalgebra structure and $\AA$ governs the algebra structure. Under some hypothesis a $\CCC^c\textrm{-}\AA$-bialgebra $\HH$ is cofree and its primitive part is governed by an operad $\PP$. We assume that the triple of operads $(\CCC, \AA, \PP)$ is good (in particular $\AA \cong \CCC \circ \PP$ as \sms). As a consequence there is an isomorphism $\HH \cong \CCC^c(\Prim \HH)$. Once such an isomorphism is chosen, we call this data a \emph{combinatorial $\CCC^c\textrm{-}\AA$-bialgebra}. In this framework one can generalize the results of this paper as follows. There exists an operad $\QQQ$ and a morphism of operads $\PP\to \QQQ$ such that the $\PP$-algebra of $\Prim \HH$ can be lifted to a $\QQQ$-algebra structure.


\end{document}